\documentclass{article}
\usepackage{amsfonts,amssymb,latexsym,amsmath}
\newcounter{num}[section]
\newcommand{\Num}{\refstepcounter{num}%
\textbf{\arabic{section}.\arabic{num}}}

\newcommand{\Theorem}{\textbf{Theorem~}}
\newcommand{\Proof}{\textbf{Proof}}
\newcommand{\Def}{\textbf{Definition~}}

\newcommand{\Lemma}{ \textbf{Lemma~}}
\newcommand{\Ex}{  \textbf{Example~}}
\newcommand{\Remark}{\textbf{Remark}}
\newcommand{\Prop}{\textbf{Proposition~}}
\newcommand{\Cor}{ \textbf{ Corollary~}}

\newcommand{\Ax}{{\mathfrak A}}
\newcommand{\Bx}{{\mathfrak B}}
\newcommand{\tx}{{\mathfrak t}}
\newcommand{\nx}{{\mathfrak n}}
\newcommand{\tCx}{\widetilde{{\mathfrak C}}}
\newcommand{\Cx}{{\mathfrak C}}

\newcommand{\Kc}{{\cal K}}
\newcommand{\Oc}{{\cal O}}
\newcommand{\Xc}{{\cal X}}
\newcommand{\Fc}{{\cal F}}

\newcommand{\Pc}{{\cal P}}
\newcommand{\Qc}{{\cal Q}}
\newcommand{\Rc}{{\cal R}}
\newcommand{\Sc}{{\cal S}}
\newcommand{\Yc}{{\cal Y}}
\newcommand{\Lc}{{\cal L}}

\newcommand{\Ch}{{{\mathfrak C}{\mathfrak h}}}
\newcommand{\al}{{\alpha}}

\newcommand{\la}{{\lambda}}

\newcommand{\tG}{{\widetilde{G}}}
\newcommand{\spann}{{\mathrm{span}}}
\newcommand{\supp}{{\mathrm{supp}}}

\newcommand{\eps}{{\varepsilon}}

\newcommand{\UTn}{{\mathrm{UT}}(n,\Fq)}
\newcommand{\utn}{{{\frak u}{\frak t}(n,\Fq)}}
\newcommand{\ut}{{\frak u}{\frak t}}
\newcommand{\UT}{{\mathrm{UT}}}
\newcommand{\Bn}{{\mathrm{B}_n}}
\newcommand{\sB}{{\mathrm{B}}}
\newcommand{\TnF}{{\mathrm{T}(n,\Fq)}}

\newcommand{\Irr}{{\mathrm{Irr}}}
\newcommand{\Ad}{{\mathrm{Ad}}}
\newcommand{\SCU}{{\mathrm{SCU}}}
\newcommand{\SCB}{{\mathrm{SCB}}}
\newcommand{\Res}{{\mathrm{Res}}}
\newcommand{\Ind}{{\mathrm{Ind}}}
\newcommand{\SInd}{{\mathrm{SInd}}}
\newcommand{\rt}{{\mathrm{right}}}

\newcommand{\NS}{{\mathrm{NS}(X)}}
\newcommand{\NSn}{{\mathrm{NS}_n(X)}}

\newcommand{\NPS}{{\mathrm{NPS}(X,Y)}}
\newcommand{\NPSn}{{\mathrm{NPS}_n(X,Y)}}
\newcommand{\Defl}{{\mathrm{Def}}}
\newcommand{\Cr}{{\mathrm{Cr}}}

\newcommand{\row}{{\mathrm{row}}}
\newcommand{\col}{{\mathrm{col}}}
\newcommand{\st}{{\mathrm{st}}}
\newcommand{\Inf}{{\mathrm{Inf}}}
\newcommand{\Fq}{{\Bbb F}_q}

\newcommand{\Cb}{{\Bbb C}}
\newcommand{\Zb}{{\Bbb Z}}
\newcommand{\Qb}{{\Bbb Q}}

\renewcommand{\leq}{\leqslant}
\renewcommand{\geq}{\geqslant}

\begin{document}
\Large

\title{Supercharacters of unipotent and solvable  groups }
\author{A.N.Panov
\thanks{The reseach is supported by the grant RSF-DFG 16-41-1013}}
\date{}

 \maketitle

\begin{abstract}
The notion of the supercharacter theory was introduced by P.Diaconis and I.M.Isaaks in 2008.  In this paper we review the main statements of the general theory, we observe the construction of supercharacter theory for algebra groups and the theory of basic characters for the unitriangular groups over the finite field.   Basing on the previous  papers  of the author, we
construct the supercharacter theory for the finite groups of  triangular type. We characterize   the structure of  Hopf algebra of supercharacters for the triangular group over the finite field.
\end{abstract}

\section{Introduction}

Traditionally, the main problem of the representation theory of finite groups is as follows:
given a finite group, to classify all its irreducible representations (characters). It turns out that for some finite groups this problem is extremely difficult, its solution is unknown and there are no
ideas how to solve it in future.    The main example of these groups is the unitriangular group  $\UT_n=\UT(n,\Fq)$ over the finite field of  $q$ elements. This group is an unipotent matrix group, and it is known   the orbit method of A.A.Kirillov is valid for these groups \cite{Kir1,Kir2,Ka}. According to this method, there exists one to one correspondence between the irreducible representations and the  orbits of  coadjoint representation.  This correspondence plays an important role in the representation theory because it enables
 to  produce a solution of the decomposition problem for the representations obtained by restriction and induction  from the irreducible ones. Moreover, the character of irreducible representation can be directly presented as a  sum taken over elements of  the coadjoint orbit (the formula of A.A.Kirillov). However,   it is not easier to classify the coadjoint orbits than to classify the  irreducible representations.  It is reasonable to say that the orbit method provides the equivalence  of categories, but not the classification of  irreducible representations.

In the series of papers   \cite{A1,A2,A3,A4}, in 1995-2003, C.Andr\'{e} investigated the theory of basic characters of the unitriangular group $\UT_n$. Although  the basic characters are not irreducible in general, the system of basic characters has many common features with the system of irreducible characters:  these characters are pairwise disjoint, they are constant on some system of subsets (called the system of basic subvarieties in $\UT_n$); the number of basic characters is equal to the number of basic subvarieties, this enables to construct the quadratic table of values (the  basic  character table).  Each basic subvariety in  $\UT_n$ is a union of the classes of conjugate elements. One can correspond the system of basic characters to some partition of the dual space   $\ut^*_n$ into the subsets, which C.Andr\'{e} called the basic subvarieties in $\ut^*_n$ (the analog of the orbit method). In contrast to
 irreducible representations, the basic characters admit the exact presentation in terms of the basic subsets in the system of positive roots of  series  $A_n$.
Observe that C.Andr\'{e} defined a basic subvariety in  $\ut^*_n$  as a sum of elementary coadjoint orbits. Apparently, this approach is universal and is valid for the other series of simple Lie algebras, but it is difficult for calculating. In  his thesis, Ning Yan  \cite{Yan} observed that for the series $A_n$ the basic subvarieties in $\ut^*_n$ coincide with the orbits of double left-right action of  the group  $\UT_n$ on $\ut_n^*$.

In the paper  ~\cite{DI}, in 2008, P.Diaconis and I.M.Isaaks presented the notion of a supercharacter theory of an arbitrary finite group. They set the  properties of the basic characters of the unitriangular group  as  axioms of a supercharacter theory.
Roughly speaking, to construct a supercharacter theory for a given finite group is to construct the system of disjoint characters (supercharacters)  $\chi_1,\ldots,\chi_m$ and the partition of the group into subsets (superclasses)  ~$K_1,\ldots,K_m$ such that the supercharacters are constant on  superclasses.     The main goal is to construct the supercharacter theory which produces the best approximation of the theory of irreducible representations.  In the paper  \cite{DI}, ~P.Diaconis and I.M.Isaaks constructed the supercharacter theory for algebra groups, its partial case is the theory of basic characters for $\UT_n$.

 Many papers were devoted to different supercharacter theories. Observe a few of them:
  supercharacters for abelian groups and their application in the number theory  \cite{Number-1, Number-2},
the superinduction for the algebra groups \cite{ThiemB, ThiemRest, MT}, the suparcharacter theory for Sylov subgroups  in orthogonal and symplectic groups over a finite field \cite{AN}, application in the random walk problem on groups \cite{Walk}, the supercharacter theory  for semidirect products \cite{H}, characterization of the Hopf algebra of supercharacters for the unitriangular group  \cite{VERY,Ben}. One can found the bibliography in the paper  \cite{VERY}.

The most of papers are devoted either to general questions of supercharacter theories, or to constructions of  supercharacter theories for the groups closely related to  $\UT_n$.  One of the main open problems is to enlarge the list, to construct an appropriate  supercharacter theory for such groups as  the parabilic subgroups in finite  Chevalley groups.
Observe that the Mackey method of classification  of irreducible representations of semidirect products apparently failed for supercharactes. This is the main  obstacle here. The general approach of the paper
\cite{H} provides rather rough supercharacter theory, which can be improved in examples.

The most of this paper is a survey.
Based on the paper  ~\cite{DI} of P.Diaconis and I.M.Isaaks we observe the general questions of supercharacter theories and construct the supercharacter theory of algebra groups (see section \ref{general} and \ref{secuni}). We also study the theory of basic characters of   C.Andr\'{e} and characterize the Hopf algebra of supercharacters of the unitriangular group following the paper  \cite{VERY} (see sections \ref{secuni} and \ref{huni}).
In section \ref{supersolv}, we  expound the supercharacter theory for the finite groups of triangular type following the author papers  \cite{P1,P2,P3}.   In section  \ref{superhopftr}, we characterize the Hopf algebra of supercharacters of the triangular group in terms of the Hopf algebra of partially symmetric functions in noncommuting variables.

\section{Supercharacters and superclasses}\label{general}
  Let $G$ be a finite group, ~$1\in G$ be its identity element, ~ $\Irr(G)$ be the set of all  irreducible characters (representations) of the group $G$.
   Given two partitions  ~
 \begin{equation}\label{partitionIrrG}
 \Irr(G) = X_1\cup \cdots \cup X_m,\quad X_i\cap X_j=\varnothing,
 \end{equation}
  \begin{equation}\label{partitionG}
G = K_1\cup \cdots \cup K_m,\quad K_i\cap K_j=\varnothing.
\end{equation}
 Observe that the partitions have equal number of components.
 We correspond $X_i$ to  the character of the group  $G$ by the formula
 \begin{equation}\label{sigma}
   \sigma_{i} = \sum_{\psi\in X_i} \psi(1) \psi.
 \end{equation}
 \Def\Num\label{twopart}. Two partitions  $\Xc=\{X_i\}$ and $\Kc=\{K_j\}$ is said to define a supercharacter theory of the group $G$ if each character  $\sigma_i$ is constant on each $K_j$. In this case,   $\{ \sigma_i\}$  are called \emph{supercharacters}, and $\{K_j\}$  \emph{superclasses}.  The table of values  $\{ \sigma_i(K_j)\}$ is called the supercharacter table.

 The next proposition  is suitable for a construction of examples  of supercharacter theories.\\
\Prop\Num\label{dilemma}~ \cite[Лемма 2.1]{DI}. Let we have the system of disjoint characters   $\Ch=\{\chi_1,\ldots,\chi_m\}$ and
  the partition $\Kc=\{K_1,\ldots, K_m\}$ of the group $G$. Suppose that each character  $\chi_i$ is constant on each  $K_j$.  Let  $X_i$ denote the \emph{support of character} $\chi_i$ (i.e. the set of all irreducible components of  $\chi_i$).  Then the following conditions are equivalent:\\
1) $\{1\}\in \Kc$,\\
2)  the system of subsets   $\Xc=\{X_i\}$ is a partition of $\Irr(G)$; the partitions  $\Xc$ and $\Kc$ define  a supercharacter theory of the group  $G$. Moreover, each $\chi_i$ equals to  $\sigma_i$ up to a constant multiplier.
Simplifying language, we  refer to  $\chi_i$ as a \emph{supercharacter}.\\
\Remark. Recall that the characters  are disjoint  (i.e. their supports  do not intersect) whenever they are orthogonal.\\
\Proof.
Assume that the condition  1) is fulfilled. Any system of disjoint characters is linearly independent. The number of characters from  $\Ch$ equals to the number of sets from  $\Kc$. Hence,    $\Ch$ is a basis in the space of all complex valued function on  $G$ constant on subsets from  $\Kc$. The regular character  $\rho(g)$ equals to  $|G|$ if $g=1$ and  equals to  $0$ if $g\ne 1$. Since  $\{1\}\in \Kc$, the character   $\rho(g)$ is constant on the subsets from  $\Kc$. We obtain
$$\rho = \sum_{i=1}^m a_i\chi_i,~~ \mbox{где}~~ a_i\in \Cb.$$ On the other hand, any irreducible character  occurs in decomposition of  $\rho$ with multiplicity equals to its degree. Therefore, any irreducible character occurs in decomposition of exactly one  $\chi_i$ and $a_i\chi_i = \sigma_i$, ~ $a_i\in\Qb^*$. This proves  2).

Assume that the condition 2) is fulfilled. The identity element of the group belongs to exactly one of the subsets of $\Kc$, say  $K_1$. Since the characters  $\chi_i$ are constant on superclasses,   $\chi_i(g)=\chi_i(1)$ for each  $g\in K_1$ and $1\leq i\leq m$. The condition  2)  implies $\rho(g) = \rho(1)$ and then  $g=1$ and $K_1=\{1\}$. ~$\Box$ \\
\Cor\Num. The system of supercharacters  $\Ch=\{\chi_i\}$ is  determined up to constants multipliers by the partition  $\Kc=\{X_i\}$.

Present some examples of supercharacter theories.\\
\Ex\Num. The system of irreducible characters  $\{\chi_i\}$ and the system of classes of conjugate elements
$\{K_i\}$.\\
\Ex\Num. Two supercharacters  $\chi_1=1_G$, ~ $\chi_2=\rho-1_G$ (here $\rho$ is the character of regular representation) and two superclasses $K_1=\{1\}$, ~ $K_2=G\setminus \{1\}$. Here is the table of supercharacters
\begin{center}\begin{tabular}{|c|c|c|}
                         \hline
                          & $\chi_1$ &  $\chi_2$\\
                          \hline
      $K_1$ & 1  &  $|G|-1$\\
      \hline
                $K_2$ & 1  &  -1\\
              \hline
               \end{tabular}
\end{center}
\Ex\Num \label{examthree}. ~ $G=C_4$ is the cyclic group of 4th order. Below we present the table of irreducible characters and the one of supercharacter tables.
\begin{center}
\begin{tabular}{|c|c|c|c|c|}
                         \hline
                          & $\chi_0$ &  $\chi_1$ &$\chi_2$ & $\chi_3$\\
                          \hline
      1 & 1  &  1  &1  & 1  \\
      \hline
                $g $& 1  &  i &-1&-i\\
              \hline
              $g^2$ & 1  &  -1 &1&-1\\
              \hline
              $g^3$ & 1  &  -i &-1& i\\
              \hline
                             \end{tabular}
               \hspace{3cm}\begin{tabular}{|c|c|c|c|}
                         \hline
                          & $\chi_0$ &  $\chi_1+\chi_3$& $\chi_2$\\
                          \hline
      1 & 1  &  2&1\\
      \hline
                $\{g, g^3\}$ & 1 &  0& -1\\
              \hline
              $g^2$ & 1 &  -2& 1\\
              \hline
                             \end{tabular}
\end{center}
\Ex\Num\label{examfour}. If the group  $\Gamma$ acts on  $G$ and $\Irr(G)$, and these actions are agree  $\chi^a(g^a)=\chi(g)$, then two partitions
$$\chi^\Gamma=\sum_{a\in \Gamma}\chi^a, \quad\quad K^\Gamma=\bigcup_{a\in \Gamma} K^a,$$ where $\chi$ (respectively, $K$) run through the set of representatives of  $\Gamma$-orbits in  $\Irr(G)$ (respectively, the classes of conjugate elements in $G$), give rise to  a  supercharacter theory of the group  $G$. The partial case in presented in the previous example: the group  $\Gamma=\Zb_2$; its generator acts on  $G$  and  $\Irr(G)$ by raising to the power three.

For the groups of small order, it is possible to enumerate  all supercharacetr theories.
It was proved in  the paper  \cite{Burkett} that there exist only three groups with exactly two supercharacter theories: the cyclic group  $Z_3$,  the symmetric group  $S_3$ and the simple group  $\mathrm{Sp}(6,2)$. For different supercharacter theories see the paper  \cite{Ash}.

Turn to formulating some general statements on supercharacters and superclasses.\\
\Lemma\Num\label{firstlemma}. If $\{X_1,\ldots , X_m\}$ is a partition of $\Irr(G)$, and $\{ K_1,\ldots, K_t\}$ is a partition of  $G$, and  $\sigma_i$, ~$1\leq i\leq m$, is constant on each  $K_j$, ~$1\leq j\leq t$, then $m\leq t$.
 \\
 \Proof. The system of characters $\{\sigma_i\}$ is linearly independent and is contained in the $t$-dimensional space of all complex valued functions on  $G$ constant on   $\{K_j\}$.  $\Box$

For any subset $K\subseteq G$,  we denote
$$\widehat{K} =\sum_{g\in K} g.$$
The center $Z(\Cb G)$  of the group algebra $\Cb G$ is a direct sum
$$Z(\Cb G)  = \bigoplus_{\psi\in\Irr(G)} \Cb e_\psi,$$
where $\{e_\psi\}$ is a system of primitive idempotents.
The  primitive idempotent   $e_\psi$, associated with the irreducible character  $\psi$, can be calculated by the formula $$
e_\psi = \frac{\psi(1)}{|G|}\sum_{g\in G} \overline{\psi(g)}g.
$$
We correspond  $X\subseteq \Irr(G)$ to the central idempotent
\begin{equation}\label{fff}
f_X = \sum_{\psi\in X} e_\psi.\end{equation}
\Prop\Num\label{hatk}~\cite{DI}. Let the partitions  $\Xc$  and  $\Kc$ define a supercharacter theory, and $f_i=f_{X_i}$.
Then  the system  $\{ \widehat{K_j}:~ 1\leq j\leq m\}$ forms a basis in the algebra $$A(\Xc,\Kc)=\spann\{f_i:~ 1\leq i\leq m\}.$$
\Proof.   Direct calculations lead to the equality
\begin{multline}\label{fpsi}
 f_i = \sum_{\psi\in X_i} e_\psi= \sum_{\psi\in X_i}\frac{\psi(1)}{|G|}\sum_{g\in G} \overline{\psi(g)}g\\
 = \frac{1}{|G|}\sum_{g\in G}\left( \sum_{\psi\in X_i} \psi(1)\overline{\psi(g)}\right) g =  \frac{1}{|G|}\sum_{j=1}^m \overline{\sigma_i(K_j)}\,\widehat{K_j}.
\end{multline}
Since the systems  $\{f_i\}$ and $\{\widehat{K_j}\}$ are of equal cardinality and linearly independent, they generate the common subspace.  $\Box$\\
\Cor\Num. If the partitions  $\Xc$  and  $\Kc$ define a supercharacter theory, then each  $K_j$  in invariant with respect to conjugation (i.e. it can be decomposed into union of classes of conjugate elements). \\
\Proof. The Proposition  \ref{hatk} implies the elements  $\{\widehat{K_j}\}$ are contain in the center of the group algebra  $\Cb G$. $\Box$  \\
\Prop\Num\label{common} ~\cite{DI}. In any supercharacter theory
the partitions $\Xc$  and $\Kc$ are uniquely determine each other.\\
\Proof.
1)  The partition  (\ref{partitionIrrG}) of the set  $\Irr(G)$ defines an equivalence relation on the group  $G$ such that  $x\sim y$ if $\sigma_i(x)=\sigma_i(y)$ for any $1\leq i\leq m$. This relation produces the partition  $\Kc^0$ of  group $G$.  Recall that   $\{\sigma_i\}$ is a basis in the space of complex valued functions constant on the partition $\Kc$. The partition  $\Kc$ is sharper than  $\Kc^0$. Therefore  $m=|\Kc|\geq |\Kc^0|$. On the other hand,  according to Lemma  \ref{firstlemma} $ |\Kc^0|\geq |\Xc|=m$. Then  $|\Kc^0|=m$  and $\Kc^0=\Kc$.\\
2) The partition  (\ref{partitionG}) of the group  $G$ defines the subalgebra  $\spann\{\widehat{K_j}\} = \spann\{f_i\}$.
Its basis of orthogonal idempotents is uniquely determined.  Therefore, the partition  $\Xc$ is uniquely defined by  $\Kc$. $\Box$\\
\Prop\Num ~\cite{DI}. The principal character  $1_G$ is a supercharacter  (recall that the character is principal if it is identically equal to one).\\
\Proof. Assume the opposite.  Let $1_G\in X_1$ and $X_1\ne\{1_G\}$. Then the partition  $\Xc^0$ that differs from $\Xc$ by decomposition  $$X_1=\{1_G\}\cup \left(X_1\setminus \{1_G\}\right),$$
forms  the pair $(\Xc^0,\Kc)$ obeying conditions of Proposition  \ref{dilemma}.
 However, $|\Xc^0|>|\Kc|$. This contradicts the conclusion of Lemma \ref{firstlemma}. $\Box$\\
\Prop\Num ~\cite{DI}. Let $\tau\in\mathrm{Aut}(\Cb)$. Then, for  any supercharacter theory  $(\Xc,\Kc)$, the automorphism  $\tau$ acts as a permutation of the partition  $\Xc$.\\
\Proof. The system $\Xc^\tau=\{ X_i^\tau\}$ is a partition of  $\Irr(G)$. The pairs $(\Xc,\Kc)$ and  $(\Xc^\tau,\Kc)$ defines two supercharacter theories with common  $\Kc$. By Proposition  \ref{common}, $\Xc=\Xc^\tau$. ~$\Box$\\
\Prop\Num ~\cite{DI}. Let  $(r, |G|)=1$. Then the bijection  $g\to g^r$ permutates superclasses.\\
\Proof. The bijection  $g\to g^r$ transform the partition  $\Kc$ into new partition  $\Kc^r=\{K_i^r\}$.
Consider the primitive root of unity  $\epsilon$ of degree $|G|$ and the automorphism   $\tau$ of the field $\Cb$ such that
  $\tau(\epsilon)=\epsilon^r$. The previous proposition implies that  $\tau$ acts on the partition  $\Xc$ by permutation and
 $\sigma_i(K_j^r)=\sigma_i^\tau(K_j)$. The pairs  $(\Xc,\Kc)$ and  $(\Xc,\Kc^r)$
define two supercharacter theories with common  $\Xc$. By Proposition  \ref{common}, $\Kc=\Kc^r$. ~$\Box$\\
\Cor\Num\label{inv}. The partition  $\Kc$ coincides with the partition  $\Kc^{-1}=\{K_j^{-1}\}$.\\
\Def\Num. Let  $A$ be a subalgebra in  $\Cb G$. The subalgebra  $A$ is called a Schur subalgebra  if there exists the partition
$\Kc=\{K_i\}$ of the group  $G$ such that \\
1)~ $\widehat{\Kc}=\{\widehat{K}_i\}$ is a basis in  $A$,~~ 2)~ $ \{1\}\in\Kc$,~~
3)~ $\Kc^{-1}=\Kc$.\\
\Theorem\Num~ \cite{H}. The map   $(\Xc,\Kc)\to A(\Xc,\Kc)$  establishes one to one correspondence between the set of supercharacters and the set  of Schur subalgebras lying in the center  $Z(\Cb G)$ of the group algebra.\\
\Proof. From Proposition \ref{hatk} and Corollary  \ref{inv} we conclude that  $A(\Xc,\Kc)$ is a Schur subalgebra lying in the center $Z(\Cb G)$. Let us show the opposite statement. Let  $A$ Schur subalgebra in $Z(\Cb G)$. The algebra  $Z(\Cb G)$ is a direct sum of a few copies of the field  $\Cb$. The subalgebra  $A$ contains 1 and is also a direct sum of a few copies of  $\Cb$. There exists a partition $\Xc=\{X_i\}$ of the set   $\Irr(G)$  such that the system of idempotents $\{f_i\}$ (see (\ref{fff}))  is a basis of $A$. The bases   $\{f_i\}$  and $\{\widehat{K}_i\}$ have equal cardinality. Then  $|\Xc|=|\Kc|$. By the formulas  (\ref{fpsi}, we obtain  $\sigma_i(K_j)=\mathrm{const}$. The partitions  $(\Xc,\Kc)$ define a supercharacter theory for the group  $G$, and $A=A(\Xc,\Kc)$. ~$\Box$

 Observe that we don't  apply the condition 3) of the definition of Schur algebra in the secod part of the proof.   Therefore, any subalgebra contained in $Z(\Cb G)$ and satisfying conditions 1) and 2)
 is a Schur algebra.

 The cited below theorem is an analog of the statement for irreducible characters.\\
\Theorem\Num~ \cite{DI}. Let $\chi$ be a supercharacter, and $K$ be a superclass. Then  $\frac{\chi(g)|K|}{\chi(1)}$ is an integer algebraic number.

Let us expound the construction of  $*$-product  of supercharacters from the paper  ~\cite{H}, which plays an impotent role in constructions of different supercharacter theories.  Let $G$ be a finite group,  $N$ be a normal subgroup in $G$, ~$\pi:G\to G/N$ be  the natural projection.

A supercharacter theory $(\Xc,\Kc)$ of $N$ is called \emph{$G$-invariant} if the action  $G$ on $N$ and   $\Irr(N)$ preserve the partitions $\Xc$  and $\Kc$.
Observe that, in a given supercharacter theory, the partitions  $G$ and $\Irr(G)$ are uniquely determine each other  (see Proposition \ref{common}). Therefore, it is sufficient to require invariance of the only one of two partitions  ($\Xc$ or $\Kc$).

Let we have the  $G$-invariant supercharacter theory  $(\Xc,\Kc)$ of the group  $N$, where $\Irr(N)= X_1\cup\ldots\cup \ X_m$ with
$X_1=\{1_N\}$ and
$N =K_1\cup\ldots\cup K_m$.  Let we have an arbitrary supercharacter theory  $(\Yc,\Lc)$ of the factor group  $G/N$, where
  $\Irr(G/N)=Y_1\cup\ldots\cup Y_k$ and  $G/N=L_1\cup\ldots\cup L_k$ with  $L_1=\{1\}$.
For any subset  $X\subset\Irr(N)$, we denote by  $\Ind(X,G)$ the sum of induced characters  $\Ind(\psi,G)$, where
$\psi\in X$.
Consider the partitions
\begin{equation}\label{sirrg}
  \Irr(G) = \supp(\Ind(X_2,G))\cup \ldots \cup \supp(\Ind(X_m,G))\cup Y_1\cup\ldots\cup Y_k.
\end{equation}
\begin{equation}\label{sggg}
  G= K_1\cup\ldots \cup K_m\cup\pi^{-1}(L_2)\cup\ldots \cup \pi^{-1}(L_k).
 \end{equation}
 \Theorem\Num~\cite{H}. The partitions (\ref{sirrg}) and  (\ref{sggg}) give rise to a supercharacter theory of the group   $G$.

 This supercharacter theory is called a  \emph{$*$-product }   of the supercharacter theories $(\Xc,\Kc)$ of $N$ and $(\Yc,\Lc)$ of the factor group  $G/N$.
 Observe that this supercharacter theory is rather rough in examples. In some papers, the construction of $*$-product was modified to obtain a more sharp  supercharacter theory.

\section{Supercharacters of finite unipotent groups }\label{secuni}

\subsection{Theory of supercharacters for algebra groups} \label{algebragroup}

 The supercharacter theory for algebra groups was constructed by P.Diaconis and I.M.Isaaks in the paper \cite{DI}. By definition, an \emph{algebra group} is a group of the form  $G=1+J$, where $J$ is an associative  finite dimensional nilpotent algebra over the finite field  $\Fq$.    The superclass of the element  $1+x$ is defined as   $1+\omega$, where $\omega $ is a left-right  $G\times G$-orbit of the element $x\in J$.
 One can define also the left and right actions of the group  $G$ on the dual space  $J^*$ by the formulas $\la g(x) = \la(gx)$ and $g\la(x)=\la(xg)$.
Let  $G_{\la,\rt}$ be the stabilizer of   $\la\in J^*$ with respect to the right action of  $G$ on $J^*$.

Fix the nontrivial character  $t\to \eps^t$ of the additive group of  field  $\Fq$ into the multiplicative group  $\Cb^*$. The function  $\xi :G_{\la,\rt}\to \Cb^*$, defined by  $$\xi_\la(g) = \eps^{\la(g-1)},$$
is a linear character  (one dimensional representation)  of the group  $G_{\la,\rt}$. Indeed, if $g_1=1+x_1$ and $g_2=1+x_2$ are two elements in $G_{\la, \rt}$, then
$\la(x_1x_2)=0$ и
\begin{multline*}
 \xi_\la(g_1g_2) = \xi_\la(1+x_1+x_2+x_1x_2) = \eps^{\la(x_1+x_2+x_1x_2)}=\\ \eps^{\la(x_1)}\eps^{\la(x_2)}\eps^{\la(x_1x_2)}=
\eps^{\la(x_1)}\eps^{\la(x_2)}=\xi_\la(g_1)\xi_\la(g_2).
\end{multline*}
A \emph{supercharacter } of the algebra group $G$ is the induced character
\begin{equation}\label{chiag}
 \chi_\la = \Ind(\xi_\la, G_{\la,\rt}, G).
\end{equation}
  \Theorem\Num\label{agmain}~ \cite{DI}. The systems of supercharacters $\{\chi_\la\}$ and  superclasses $\{1+GxG\}$, where $\la$  and $x$  run through the systems of  representatives of  $G\times G$-orbits in  $J^*$ and $J$ respectively, give rise to  a supercharacter theory  of the group  $G$.

  Thus, the characters of the system  $\{\chi_\la\}$ are pairwise disjoint. The classification problem of irreducible representation reduces to decomposition of supercharacters into a sum of irreducible components. This problem remains open up today. It could be interesting the following result.\\
\Theorem\Num\label{intersec} ~\cite{Asite}. $(\chi_\la,\chi_\la)= |G\la\cap \la G|$.\\
\Proof.
One can realize the representation with  character $\chi_\la$ in the space  $V_\la=<f_\mu(1+x) =\eps^{ \mu(x)}: ~\mu\in G\la>$ by the formula $T_gf(s)=f(sg)$. That is, the representation with  character  $\chi_\la$  is a subrepresentation of the right regular representation of the group  $G$ in the space of complex valued functions  $\Cb[G]$.
 The representation in the space  $\Cb[G]$ is completely reducible;  $\Cb[G]$ is decomposed into a direct sum  of $V_\la$ and its complement $V_\la'$. Any  $\phi \in \mathrm{Hom}_G(V_\la,V_\la)$ extends to the intertwining operator  $\tilde{\phi}$ in the space  $\Cb[G]$ that is equal to zero on  $V_\la'$. Then  $\tilde{\phi}$ is a linear combination  of the left translation operators  $L_gf(s)=f(gs)$ on the group $G$. The operator of  left translation by  $a\in G$ preserves the subspace  $V_\la$ whenever  $\la a\in G\la$.  Then  $|\mathrm{Hom}_G(V_\la,V_\la)|=|G\la\cap \la G|$. $\Box$

Let us recall the A.A.Kirillov formula for the irreducible characters.
For each   $\la\in J^*$ one can construct the irreducible character  $\Psi_{\la}$  of the algebra group  $G=1+J$ (see \cite{Kir1}); it can be presented in the form:
\begin{equation}\label{kkk}
  \Psi_\la (1+x) = \frac{1}{\sqrt{|\Omega^*|}}\sum_{\mu \in\Omega^*} \varepsilon^{\mu(x)},
\end{equation}
where $\Omega^*$ is the coadjoint orbit of $\la\in J^*$.
 Observe that  this formula is valid if the characteristic of the field is sufficiently great  \cite{Kir1}, \cite{Ka}.

        The functions defined by (\ref{kkk})  are called \emph{Kirillov functions}.
The number of Kirillov functions equals to the number of irreducible representations (by ~\cite[Lemma 4.1]{DI}), however,  this sets do not coincide in general.
The simplest example: the group   $\UTn$, ~ $n\ge 13$, for   $\mathrm{char}(\Fq)=2$
 ~\cite{IsKar}.

In the paper \cite{DI}, the supercharacter analog of the formula (\ref{kkk}) for the  algebra groups  was presented:
\begin{equation}\label{kirdi}
  \chi_\la (1+x) = \frac{1}{n(\la)}\sum_{\mu \in G\la G} \varepsilon^{\mu(x)},
\end{equation}
where $ \chi_\la$  is the supercharacter associated with  $\la\in J^*$,
 and $n(\la)$ is the number of right  $G$-orbits in $G\la G$. This formula is valid over a field of an arbitrary characteristic.

  The constant $n(\la)$ is of interest for the other reason.
 Recall that there is the standard supercharacter $\sigma_\la$ that is  the sum of all irreducible constituents of $\chi_\la$ with multipliers equal to their degrees (see Proposition
 \ref{dilemma}).  The supercharacter $\chi_\la$ differs from  $\sigma_\la$ by a constant multiplier. It turns out that this constant is equal to
  $n(\la)$, more precisely  $n(\la)\chi_\la=\sigma_\la$ ~\cite{DI}.

In the same paper  \cite{DI}, the other version of the formula (\ref{kirdi}) was presented in the form
\begin{equation}\label{kirdisec}
  \chi_\la (1+x) = \frac{|\la G|}{|GxG|}\sum_{y\in GxG} \varepsilon^{\la(y)}.
\end{equation}

Turn to the questions on restriction and induction in the supercharacter theory.
Let $G =1+J$ be an algebra group. If $J'$ is an arbitrary its subalgebra, then the subgroup $G'=1+J'$ as called an \emph{algebra subgroup} in  $G$. \\
It is proved in the paper \cite{DI} that the restriction of the supercharacter  $\chi_\la$ on the algebra subgroup  $G'$ is a sum of supercharacters of the subgroup  $G'$ with non-negative integer coefficients.

  Let $\phi$ be a superclass function on $G'$ (i.e.,  it is a function constant on superclasses in $G'$).  Extend  $\phi$  to a function  $\dot{\phi}$ on $G$ taking it equal to zero outside of $G'$. \\
 By definition, a \emph{superinduction }  for  $\phi$ is a function    $\SInd\,\phi$ on the group   $G$ defined by the formula
$$ \SInd\,\phi(1+x) = \frac{1}{|G|\cdot|G'|}\sum_{a,b\in G}\dot{\phi}(1+axb).$$
Easy to see that  $ \SInd\,\phi(1+x)$ is a superclass function on  $G$.
The scalar product on  $G$ (and $G'$) is defined as usual
\begin{equation}\label{scalar}
  (f_1,f_2) = \frac{1}{|G|}\sum_{g\in G} f_1(g)\overline{f_2(g)}.\end{equation}

The following theorem is a supercharacter analog of the Frobenius theorem.\\
\Theorem\Num~\cite{DI}. Let $\phi$ be a superclass funstion on  $G'$, and $\psi$ be a superclass function on  $G$. Then
$(\SInd\,\phi,\psi) = (\phi, \Res\,\psi).$

\subsection{Theory of basic characters of the unitriangular group}\label{basic}

The most impotent example of an algebra group is the unitriangular group $\UT_n=\UTn$ that consists of all upper triangular matrices the entries from the finite field  $\Fq$ and with ones on the diagonal.  Its Lie algebra  $\ut_n=\utn$ consists of all strong upper triangular matrices, and, indeed, it is an associative algebra. The group   $\UT_n=1+\ut_n$ is an algebra group. Created by C.Andr\'{e}  the theory of basic characters \cite{A1,A2,A3} is a special case of the P.Dianconis and I.M.Isaaks theory of supercharacters for algebra groups  \cite{DI}.

Introduce a few definitions.
 A \emph{root} is a pair of positive integers  $(i,j)$, where $1\leq i,j\leq n$ and $i\ne j$.
 The root  $(i,j)$ is  \emph{positive} (respectively, negative) if  $i<j$ (respectively, if $i>j$).
 The number $i$ is called a \emph{number of row} of the root  $\al=(i,j)$ (denote $i=\row(\al)$).
 Respectively, $j$ is a \emph{number of column}  (denote $j=\col(\al)$).
 The Lie algebra   $\ut_n$ has the basis of matrix units  $\{E_\al:~ \al>0\}$.

 The subset of positive roots $ D$ is called a \emph{basic subset} (after
 C.Andr\`{e})  if it has  at most  one root  in each row and each column.
 The other name of $D$ is a \emph{subset of rooks arrangement type}.
 Suppose that we have a map  $\phi:D\to \Fq^*$.
 The pair  $(D,\phi)$ is called an \emph{admissible pair} if  $D$  is a basic subset.

 For each admissible pair   $(D,\phi)$,  we   consider the element
 $$ X_{D,\phi}=\sum_{\al\in D} \phi(\al)E_\al$$
 in the Lie algebra  Ли $\ut_n$.
 The dual space   $\ut_n^*$ has the basis  $\{E^*_\al:~ \al>0\}$ dual to $\{E_\al:~ \al>0\}$.
 Denote by   $\la_{D,\phi}$ the  element of the dual space of the form
  $$ \la_{D,\phi}=\sum_{\al\in D} \phi(\al)E^*_\al.$$
 \Theorem\Num ~\cite{DI, Yan}. Let $G=\UT_n$. Then \\
 1) Any left-right orbit of the group  $G\times G$ in $\ut_n$   contains the unique element of the form  $X_{D,\phi}$ (denote the orbit by $\Oc(D,\phi)$).\\
  2) Any superclass of the group  $G=\UT_n$ contains the unique element of the form   $1+X_{D,\phi}$ (denote the superclass by  $K_{D,\phi}$).\\
 3) Any left-right orbit of the group $G\times G$ in $\ut_n^*$  contains the unique element of the form   $\la_{D,\phi}$ (denote the orbit by  $\Oc^*(D,\phi)$).

    The group $G$ decomposes into superclasses  $\{K_{D,\phi}\}$,  and  $\ut_n^*$ into the left-reght orbits  $\{\Oc(D,\phi)\}$. Following the scheme of subsection \ref{algebragroup}, we construct the supercharacter  $\chi_{D,\phi}$ by the pair   $(D,\phi)$. In the papers of C.Andr\`{e}, these supercharacters are called  the \emph{basic characters} and these superclasses are \emph{basic subvarieties} in $\UT_n$.
   The following theorem is a special case of Theorem  \ref{agmain}.\\
   \Theorem\Num.  The systems of basic characters  $\{\chi_{D,\phi}\}$ and basic subvarieties  $\{K_{D,\phi}\}$ give rise to a supercharacter theory on the group  $\UT_n$.

   We correspond the  basic subset $D$  to the graph with vertices in  $1,2,\ldots, n$; two vertices $i$ and $j$ are connected by a edge if the root  $(i,j)$ belongs to  $D$. For example,  for $D=\{(1,3), (3,6), (2,4), (4,5), (5,7)\} $, the graph has the form:

   
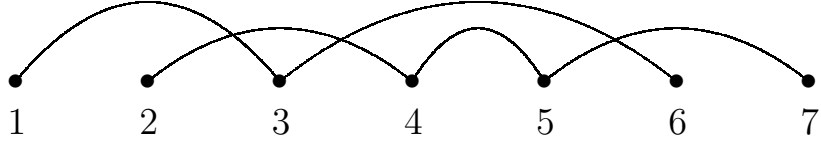
\begin{picture}(500,100)
\qbezier(50,40)(100,100)(150,40)
\qbezier(150,40)(225,100)(300,40)
\qbezier(100,40)(150,80)(200,40)
\qbezier(200,40)(225,80)(250,40)
\qbezier(250,40)(300,80)(350,40)
\put(50,40){\circle*{5}}
\put(100,40){\circle*{5}}
\put(150,40){\circle*{5}}
\put(200,40){\circle*{5}}
\put(250,40){\circle*{5}}
\put(300,40){\circle*{5}}
\put(350,40){\circle*{5}}
\put(47,20){\text{1}}
\put(97,20){\text{2}}
\put(147,20){\text{3}}
\put(197,20){\text{4}}
\put(247,20){\text{5}}
\put(297,20){\text{6}}
\put(347,20){\text{7}}

\end{picture}

We say that the roots $(i,j)$ and  $(k,l)$ \emph{form a crossing} if  $i<k<j<l$.
 Denote by  $c(D)$ the number of crossings in  the basic subset  $D$.
In the above figure,  intersections of edges match crossings of roots;  $c(D)=3$.
 For any basic subset $D$, we take
    $$ d(D)=\sum_{(i,j)\in D} (j-i-1).$$
    The Theorem  \ref{intersec} implies the following statement.\\
\Theorem\Num\label{spcros}~\cite{VERY}. $(\chi_{D,\phi},\chi_{D,\phi})=q^{c(D)}$ and $\chi_{D,\phi}(1) = q^{d(D)}$.

  If  $D$ consists of the only root  $\{\al\}$ and $c=\phi(\al)\in \Fq^*$, then we denote by
  через $\chi_{\al,c}$ its supercharacter, and by    $\Oc_{\al,c}$ and $\Oc^*_{\al,c}$ the corresponding left-right orbits in $\ut_n$ and $\ut_n^*$. Observe that, in this case,   $\chi_{\al,c}$ is an irreducible character of the group  $\UT_n$ (called an  \emph{elementary  character}),  $\Oc_{\al,c}$  is an adjoint orbit, $\Oc^*_{\al,c}$ is a coadjoint orbit (both are called  \emph{elementary orbits}).
 \\
\Theorem\Num ~\cite{A2,A3}. Let  $(D,\phi)$ be an admissible pair, ~$D=\{\al_1,\ldots,\al_k\}$ and $c_i=\phi(\al_i)$. Then\\
1) the left-right orbit $\Oc(D,\phi)$  in  $\ut_n$ is the sum of elementary orbits  $\sum_{i=1}^k \Oc_{\al_i,c_i}$;\\
2)  the left-right orbit $\Oc^*(D,\phi)$  in $\ut_n^*$ is the sum of elementary orbits $\sum_{i=1}^k \Oc^*_{\al_i,c_i}$;\\
3)  the supercharacter   $\chi_{D,\phi}$  of the group   $\UT_n$ is the product  of elementary characters  $\prod_{i=1}^k \chi_{\al_i,c_i}$.

  The problem of calculation of values of a given supercharacter on superclasses reduces to the case of elementary character. We need new definitions.

Define the operation of  partial addition in the set of positive roots by  $(i,j)=(i,l)+(l,j)$. Here, we say that  the roots $(i,l)$ and $(l,j)$ are \emph{singular}  for the root  $(i,j)$. Denote the number of singular roots for  $\al$ by  $\mathrm{Sing}(\al)$.
For any positive root  $\al=(i,j)$ and any  basic subset  $D'$,
denote by  $D'(i,j)$ the number of  $\beta\in D'$ such that  $\row(\beta)>i$ and $\col(\beta)<j$. Take  $d'(i,j)=j-i-1-|D'(i,j)|$.\\
\Theorem\Num \cite{A3, VERY}. Let $\al=(i,j)$ be a positive root, ~ $c\in\Fq^*$, and  $(D',\phi')$ is an arbitrary admissible pair. Then the value of the supercharacter $\chi_{\al,c}$ on the superclass $K_{D',\phi'}$ can be calculated  by the formula
\begin{equation}\label{elemchar}
  \chi_{\al,c}(K_{D',\phi'}) = \left\{\begin{array}{cl} q^{d'(i,j)}&, ~~\mbox{if}~~ D'\cap\mathrm{Sing}(\al)=\varnothing, ~~\al\notin D';\\
  q^{d'(i,j)}\eps^{c\phi'(\al)}&, ~~\mbox{if}~~ D'\cap\mathrm{Sing}(\al)=\varnothing, ~~\al\in D';\\
  0 &, ~~\mbox{if}~~ D'\cap\mathrm{Sing}(\al)\ne\varnothing.\end{array}\right.
\end{equation}

One can consider the basic subvarieties in  $\ut_n^*$ defined over an arbitrary field.  The questions of decomposition of basic subvarieties into the coadjoint orbits and finding generators in the field of invariants are discussed in the papers  \cite{A5,P4}.
 Further, we observe the  characterization of the basic subvarieties in terms of the tangent cones for  the Schubert varieties.

Let  $K$ be an arbitrary algebraically closed field.
The dual space  decomposes into the basic subvarieties  $\Oc^*_{D,\phi}$,
 where $D$ is a basic subset, and $\phi$ is a map $D\to K^*$.
A  \emph{basic cell}  $ V^\circ_D$ is a union
$$V^\circ_D = \bigcup_{\phi}\Oc^*_{D,\phi};$$
we refer to its  closure $V_D$ as a \emph{ basic cone}.

 Consider the flag variety  $\Fc=G/B$, where $G=\mathrm{GL}(n,K)$  and $B=T(n,K)$. It is well known that $\Fc$
is decomposed into the Schubert cells
$\Fc_w^\circ =  B\dot{w}B\bmod B,$
where $w$ is the element of the Weyl group  $W$, and $\dot{w}$ is its representative in the group $G$.
 The closure of the Schubert cell $\Fc_w^\circ$ in the Zariski topology is called a \emph{Schubert variety} $\Fc_w$.  Any Schubert variety contains the origin point  $p=B\bmod B$.

Let  $N_-$ be the group of lower triangular matrices of order $n$ with ones on the diagonal. The subset  $V= N_-B\bmod B$ is called the  \emph{main affine neighborhood} of the point  $p\in V$.
Applying the presentation   $N_-=1+\nx_-$, one can identify  the  main affine neighborhood with
$\nx_-$,  and, in its turn, $\nx_-$, using the Killing form, with $\ut_n^*$; then   $p=0$. Thus, the intersection  $\Fc_w\cap V$ is a Zariski-closed subset in $\ut_n^*$.
The tangent cone  $\mathrm{TC}_w$ of the Schubert variety  $\Fc_w$ at the point $p=0$
is a closed cone in   $\ut_n^*$ with the center at zero.\\
\Def\Num. We call the element $w\in W$  \emph{ homogeneous} if
$\Fc_w\cap V$ is a cone with the center at the point $p$ (i.e., $\mathrm{TC}_w=\Fc_w\cap V$).\\
\Theorem\Num~ \cite{P4}.  The basic cones coincide with the tangent cones of the homogeneous elements of the Weyl group.

\subsection{Irreducible constituents of supercharacters}
The classification problem of irreducible representations of the unitriangular group  reduces to the  problem of decomposition of a given supercharacter into a sum of irreducible constituents.  The last problem is also  extremely difficult in general. However, there are some successes towards a solution.

Observe that the supercharacters  $\{\chi_{D,\phi}\}$ with the common $D$ and different  $\phi:D\to \Fq^*$ are conjugate with  respect to the adjoint action of the diagonal subgroup  $H=\Fq^{*n}$. It is sufficient to obtain  a decomposition of only one supercharacter from the series, for example for   $\phi\equiv 1$.

First question, when a basic character s irreducible?
From  Theorem  \ref{spcros}, the basic character   $\chi_{D,\phi}$ is irreducible if and only if $c(D)=0$, that is, the graph of the basic set has no crossings.

The second question, when a basic character is a multiple irreducible?
By a $k$-crossing, we call the series  $i_0<i_1<i_2<\ldots<i_k<i_{k+1}<i_{k+2}$, where each pair $(i_s,i_{s+2})\in D$. The number  $k$ stands for the length of the crossing. A maximal crossing is defined in a usual way as a crossing that can't be extended. The following statement can be derived from   ~\cite[Теорема 4]{A6}.\\
\Theorem\Num. The basic character  $\chi_{D,\phi}$  is multiple irreducible if and only if each maximal crossing in $D$ has even length. In this case, its  unique irreducible component has degree  $q^e$, where  $e=d(D)-\frac{1}{2}c(D)$.

In the general case, the supercharacter $\chi_{D,\phi}$ is  a linear combination with integer non-negative coefficients of the classfunctions from the set  $\mathrm{Kir}(D,\phi)$ that consists of  the Kirillov functions associated with the coadjoint orbits from  $G\la_{D,\phi} G$ ~\cite{A7}.
 The number of  the irreducible characters in  $supp(\chi_{D,\phi})$  coincides with $|\mathrm{Kir}(D,\phi)|$~ \cite{Otto}. One can  hypothesize that the sets   $\supp(\chi_{D,\phi})$  and $\mathrm{Kir}(D,\phi)$ also coincide, but this is hardly true.

According to Higman's conjecture, the number of irreducible representations of the group  $\UT (n,\Fq)$ is a polynomial in  $q$. One can set the  similar conjectures for the number
$N_D(q)$  of irreducible constituents of the basic character $\chi_{D,\phi}$  and for the number
$N_{D,e}(q)$  of its irreducible constituents of  degree  $q^e$. This hypotheses are not proved ip today. The survey of contemporary results is contained in  ~\cite{M2}.

In the paper ~\cite{M2}, it is verified that the number   $N_{D,e}(q)$ can be calculated in terms of irreducible representations of  certain  algebra  $\tCx_D(q)$.
Denote by  $\Cr(D)$ the set of all positive roots $(i,j)$ such that there exist the roots   $(i,k)$, $(j,l)$ from $D$, where  $i<j<k<l$ (i.e., the roots   $(i,k)$ and  $(j,l)$ produce a crossing).
First, we construct the algebra  $\Cx_D(q)$ as an algebra generated by the elements  $e_{ij}$, where $(i,j)\in \Cr(D)$, with the relations
$$ e_{ij}* e_{kl} =\left\{\begin{array}{rl} e_{il},& ~~\mbox{if}~~ j=k ~~\mbox{and}~~  (i,l)\in\Cr(D),\\
0,& ~~\mbox{~~ otherwise ~~}.\end {array}\right.$$
This multiplication does not coincide with the multiplication of matrix units  $E_{ij}$ and $ E_{kl}$.
The algebra $\tCx_D(q)$ is its central extension  $\Cx_D(q)\oplus \Fq z_D$ with the multiplication
$z_D^2=0$,~ $e_{ij}z_D=z_De_{ij}=0$ and
$$ e_{ij}* e_{kl} =\left\{\begin{array}{ll} e_{il}&, ~~\mbox{if}~~ j=k ~~\mbox{and}~~  (i,l)\in\Cr(D),\\
z_D&, ~~\mbox{if}~~ j=k ~~\mbox{and}~~  (i,l)\in D,\\
0&, ~~\mbox{~~ otherwise ~~}.\end {array}\right.$$
The algebra $\tCx_D(q)$, as  $\Cx_D(q)$, is associative and nilpotent. Let  $J=\ut_n$.
Observe that
$$ \Cx_D(q)\cong {\mathfrak s}_\la/{\mathfrak l}_\la\quad \mbox{and}\quad \tCx_D(q)\cong {\mathfrak s}_\la/{\mathfrak k}_\la,$$
 where ${\mathfrak l}_\la = J_{\la,\rt}$, ~~ $ {\mathfrak s}_\la=\{x\in J: ~\la(xy)=0, ~~\mbox{ for~~ all}~~ y\in J_{\la,\rt}\}$, ~~${\mathfrak k}=J_{\la,\rt}\cap \mathrm{Ker}(\la)$.
Let  $\Irr(G,k)$ be the set of the irreducible representations of the group  $G$ of degree $k$.\\
\Theorem\Num ~\cite[Theorem 3.1]{M2}. For any basic subset  $D$ and the field  $\Fq$,  we obtain the equality
$$
N_{D,e}(q) = \frac{\sharp(\Irr(1+\tCx_D(q), q^f) - \sharp(\Irr(1+\Cx_D(q), q^f)}{q-1}.$$
where $f=c(D)-d(D)+e$. More over, if all irreducible characters of the group  $1+\tCx_D(q)$ of degree  $q^f$ are the Kirillov functions, then all irreducible constituents in $\chi_{D,\phi} $ of degree  $q^e$  are also Kirillov functions.

\section{The Hopf algebra of supercharacters for the unitriangular group}\label{huni}

\subsection{The Hopf algebra  $\NS$}\label{hau}
Let us recall the definition of the Hopf algebra of symmetric functions in non-commuting  variables ~\cite{VERY}. Let  $X=\{x_1,x_2,\ldots\}$ be a  set of non-commuting variables (alphabet). The group $S_\infty=\varinjlim S_n$ acts on the set   $X$  by  finite permutations  of the variables.  Consider the linear space  $\NS$
of all formal power series in variables of $X$ with complex coefficients and  bounded degree, and invariant with respect to the group  $S_\infty$.  The linear space  $\NS$ is a direct sum of the subspaces  $$ \NS = \sum_{n=0}^\infty \NSn $$ of all homogeneous elements of the fixed degree.  For $n=0$,   the space  $\mathrm{NS}_0(X)$ coincides with the field  $\Cb$.
The linear space $\NS$ is a  graded algebra with respect to the natural multiplication. The algebra  $\NS$ is called the  \emph{algebra of symmetric functions in non-commuting variables}.

The algebra  $\NS$ is a Hopf algebra defined as follows. Let  $X'$ and $ X''$  be two pairwise commuting copies of the alphabet  $X$. Then,  for any  $f(X)\in\NS$,
we have decomposition  $f(X'+X'')=\sum f'_k(X')f''_k(X'')$, where $X'+X''$ is a disjoint union of the alphabets  $X'$ and $X''$. By definition,  $$\Delta (f)=\sum f_k'\otimes f_k''.$$

Let  $P$ be a set partition of $[n]=\{1,\ldots, n\}$ (denote  $P\vdash[n]$).  Let $m_P$
be a sum of all monomials  $x_{i_1}x_{i_2}\cdots x_{i_n}$, where
$x_{i_s}=x_{i_t}$  whenever  $i_s$ and $i_t$ belong  to a common part of the partition  $P$. For example, of  $P= 13|2|4$, then $m_p=x_1x_2x_1x_3 + \ldots$.

The system of elements $\{m_P: ~ P\vdash [n]\}$ is a basis in   $\NSn$.
If $P\vdash [k]$ and $Q\vdash [m]$, then $m_Pm_Q\in \NSn$, ~ $n=k+m$, and
$$m_Pm_Q =  \sum m_R,$$
where the sum is taken over all $R\vdash [n]$, ~ $R\wedge ([k]|[m])=(P|Q)$.

 By a \emph{subpartition} of the partition  $P\vdash [n]$, we call the system of subsets  $P_1$ such that each its subset is one of components  of the partition $P$.
 By definition, $P=P_1+P_2$  if the partition  $P$ is a disjoint union of two subpartitions $P_1$ and $P_2$.

 For any subset  $A\subseteq [n]$, denote by  $\st $ the unique order preserving map  $\st: A\to [|A|]$. Define comultiplication as follows  $$\Delta(m_P)=\sum_{P=P_1+P_2} m_{\st(P_1)}\otimes m_{\st(P_2)}.$$
 Example. For the partition   $P=14|2|3|$ of the segment  $[4]$, we have $\Delta(m_P) = m_{14|2|3|}\otimes 1+ 2  m_{13|2}\otimes m_1 + m_{12}\otimes m_{1|2} +
m_{1|2}\otimes m_{12} +  2  m_1\otimes m_{13|2} + 1\otimes  m_{14|2|3|}$.

The linear space $\NS$ is a Hopf algebra with respect to the defined above multiplication,
comultiplication, the unit  $1\to m_\varnothing$ and the counit  $f\to f(0,0,\ldots)$.

\subsection{The Hopf algebra  $\SCU$}\label{scu}
Let  $q=2$ and  $\UT_n$ be the unitriangular group of the $n$th order defined over the field of two elements. Consider the linear space of superclass functions $\SCU_n$ on the group  $\UT_n$ and the direct sum  $$ \SCU=\sum_{n=1}^\infty \SCU_n.$$
Define the multiplication and comultiplication on  $\SCU$.  The group  $\UT_k\times \UT_m$ is a factor group of  $\UT_n$, where $n=k+m$.
Let $\phi\in \SCU_k$ and $\psi\in \SCU_m$, then, by definiiton,
  $$\phi\cdot \psi = \Inf(\phi\times\psi)\in \SCU_n,$$
  where the inflation is a composition of  $\phi\times\psi$ and the natural projection   $\UT_n\to \UT_k\times \UT_m$.

  For any partition  $T=(A_1|A_2)$ of the set  $[n]$,  there is the subgroup  $\UT_{A_1}\times \UT_{A_2}$ of the group  $\UT_n$; it consists of  all unitriangular matrices such that  $a_{ij}\ne 0$ implies  that  $i$, $j$  belong to a common component of the partition $T$.  The subgroup  $\UT_{A_1}\times \UT_{A_2}$ is natually isomorphic to  $\UT_{|A_1|}\times \UT_{|A_2|}$ (denote this isomorphism by $\pi_T$).
 Let  $\chi\in \SCB_n$. Then, by definition,
 $$\Delta(\chi) = \sum_{T=(A_1|A_2)\vdash [n]} { ^T\Res}^{\UT_n}_{\UT_{|A_1|}\times \UT_{|A_2|}}(\chi),$$
where $^T\Res(\chi)(g)=\Res(\chi)(\pi_T^{-1}(g))$.

We correspond  each partition $P\vdash [n]$ to the basic subset  $D_P$ as follows. The root $(i,j)$, where $1\leq i<j\leq n$, belongs to  $D_P$ if  $i,j$ are the adjacent elements in a common component of the partition  $P$. For example, if $P=(135|24)$, then  $D_P=\{(1,3), (3,5),( 2,4)\}$.

Since  $q=2$, we have  $\phi\equiv 1$;  each partition  $P$ corresponds to a unique superclass  $K_P$ and supercharacter  $\chi_P$. Denote by $\kappa_P$ the characteristic function on the superclass  $K_P$.
The system   $\{\kappa_P:~ P\vdash [n]\}$ is a basis in  $\SCU_n$.\\
\Theorem\Num \cite[Theoremа 3.2, Corollary 3.3]{VERY}. The linear space  $\SCU$ is a Hopf algebra. The map   $\kappa_P\to m_P$, ~ $P\vdash [n]$ can be linearly extended  to the isomorphism of the Hopf algebra  $\SCU$ onto the Hopf algebra  $\NS$.
\\
\Remark~\Num. For the arbitrary $q$, the Hopf algebra  $\SCU$ is isomorphic to the "coloured"  $\NS$ ~\cite[Theorem 3.6]{VERY}.\\
\Remark~\Num\label{hopfave}.
For the arbitrary $q$,  one can consider the supercharacter theory constructed by  averaging  with respect to the action of group of diagonal matrices  $\Gamma$
(in  spirit of example 4 of section \ref{general}). Then each  partition  $P\vdash [n]$ corresponds to the unique supercharacter   $\chi^\Gamma_P$ and superclass  $K^\Gamma_P$. The map that sends the characteristic function   $\kappa^\Gamma_P$  of the superclass  $K^\Gamma_P$ to the element  $m_P$  linearly extends to the  isomorphism  of the Hopf algebra  $\SCU$ onto the Hopt algebra  $\NS$~\cite{Asite}.

\subsection{The dual Hopf algebra  $\SCU^*$}\label{dualu}

The linear space   $\SCU^*=\sum_{n=0}^\infty \SCU^*_n$ is a Hopf algebra with respect to operations dual to the one in  $\SCU$. The scalar product  $(\cdot,\cdot)$  on the group   $G$ extends to the scalar product on  $\SCU$. Identify  $\SCU^*$ and  $\SCU$ as  linear spaces; let us calculate the dual operations.

 Let $\phi\in \SCU_k$, ~ $\psi\in \SCU_m$ and $n=k+m$. For the partition  $T=(A|A^c)$, the operation dual to  $^T\Res$  is the operation $^T\SInd: \SCU_k\otimes \SCU_m\to \SCU_n$ defined as
$$^T\SInd(\phi\times\psi)= \SInd (\phi\pi_A\times\psi\pi_A)$$
 The multiplication in $\SCU^*$ is defined by the formula
 \begin{equation}\label{ssind}
\phi\cdot\psi = \sum_{T=(A_1|A_2),~ |A_1|=k, ~|A_2|=n-k} { ^T\SInd}^{\UT_n}_{\UT_{|A_1|}\times \UT_{|A_2|}}(\phi\times\psi).
\end{equation}
By definiton, the comultiplication in  $\SCU^*$ is the deflation  $\Defl$ that is dual to inflation $\Inf$; it is defined as follows.
Let  $\tau$ be the natural projection  $\UT_n\to \UT_k\times\UT_m$ and $\chi\in \SCU_n$. Then
\begin{equation}\label{defl}
\Delta(\chi)(u,v) = \Defl(\chi) (u,v) = \frac{1}{| \tau^{-1}(1)|}\sum_{x\in\tau^{-1}(u,v)}\chi(x),
\end{equation}
where $u\in \UT_k$, ~ $v\in\UT_m$.\\
\Theorem\Num~ \cite{VERY}. The operations of superinduction  (\ref{ssind}) and deflation  (\ref{defl}) give rise to the structure of  dual Hopf algebra  $ \SCU^*$ on  the linear space $ \SCU$.

The system of elements $\{\kappa_P:~ P\vdash [n]\}$ is a basis in  $\SCU_n$.
The dual basis  $\{\kappa_P^*:~ P\vdash [n]\}$ consists of the elements    $\kappa_P^*=z_P\kappa_P$, where
$$z_P=\frac{|\UT_n|}{|\UT_nX_{D_P}\UT_n|}.$$
Let us calculate the operations of $ \SCU^*$ on the  elements of the basis.
For $P\vdash[k]$ and $Q\vdash[m]$, we have
$$\kappa_P^*\cdot\kappa_Q^* = \sum \kappa^*_{(\st_{A_1}^{-1}(P)|\st_{A_2}^{-1}(Q))},$$
where the sum is taken over all the partitions  $(A_1,A_2)\vdash [n]$ such that  $|A_1|=k$, ~$|A_2|=n-k$.

 Let  $P\vdash[n]$.  Then
  $$\Delta(\kappa_P^*) = \sum_{k=0}^n \kappa^*_{P_{[k]}}\otimes \kappa^*_{P_{[k]^c}},$$
where $P_{[k]}$ and  $P_{[k]^c}$ are the intersections of the partition  $P$ with $[k]$ and its complement  $[k]^c$ in  $[n]$.

\section{ Supercharacters of finite solvable groups}\label{supersolv}

\subsection{Supercharacter theory for finite groups of triangular type}\label{supertrtr}
In this subsection, we construct the supercharacter theory the finite groups of triangular type following the author papers  \cite{P1,P2, P3}.

Let  $H$ be a group, and  $J$ be an associative algebra defined over a field $k$. We assume that there are defined  the commuting  left   $h,x\to hx$ and right  $h,x\to xh$ linear actions of the group  $H$ on  $J$. We suppose that,  for any  $h\in H$ and $x,y\in J$, the following conditions are fulfilled:
\begin{enumerate}
\item $h(xy)=(hx)y$ и $(xy)h=x(yh)$,
\item $x(hy)=(xh)y$.
\end{enumerate}
The formula \begin{equation}\label{mult}
g_1g_2=(h_1+x_1)(h_2+x_2)=h_1h_2+ h_1x_2+x_1h_2+x_1x_2
\end{equation}
defines an associative operation on the set  $$G=H+J=\{h+x:~ h\in H, ~ x\in J\}.$$

If  $J$ is a nilpotent algebra over the field  $k$, then  $G$ is a group with respect to the operation  (\ref{mult}). If the group $H$ is finite, and   $k=\Fq$, and $J$ is a finite dimensional  nilpotent algebra over  $k$,  then  the group  $G$ is  finite. \\
\Def\Num. Under the above requirements, we call the  group  $G$ a \emph{finite group of triangular type} if the group  $H$ is abelian and  $\mathrm{char}\, k$ does not divide  $|H|$.

Let $G=H+J$ be a finite group of triangular type. The group algebra  $kH$ is commutative and semisimple, by Maschke's theorem. Therefore,  $kH$ is a sum of fields.
There exist the system of primitive idempotents  $\{e_1,\ldots, e_n\}$ such that   \begin{equation}\label{gralg} kH=k_1e_1\oplus\ldots \oplus k_ne_n, \end{equation} where $k_1,\ldots,k_n$ are  extensions of the field  $k$. Any idempotent from  $kH$ is a sum of primitive idempotents.

The direct sum  $A=kH\oplus J$ has structure of an algebra with respect to the multiplication  (\ref{mult}). The group  $G$ is a subgroup in the group  $A^*$ of invertible elements of the algebra  $A$, see the example \ref{exthree} below.
Observe that the group  $G$  is decomposed into the product  $G=HN$ of its subgroup $H$ and the normal subgroup   $N=1+J$, which is an algebra group.\\
\Ex\Num. The algebra group $G=1+J$.\\
\Ex\Num\label{extwo}. $G=\left\{\left(\begin{array}{cc}a&b\\
0&1\end{array}\right)~~ a,b\in \Fq,~ a\ne 0\right\}$.\\
\Ex\Num\label{exthree}.
Let $A$ be an associative finite dimensional algebra with unit over the finite field   $\Fq$ of $q$ elements \cite[\S 6.6]{Pi}. By definition,  an algebra $A$ is  reduced if the its factor algebra over  the radical    $J=J(A)$ is a direct sum of division algebras. According to Wedderburn's theorem  \cite[\S 13.6]{Pi}, any division algebra over a finite field is commutative. Then the algebra  $A/J$ is commutative.
There exists a semisimple subalgebra $S$ such that
$A = S\oplus J$ (see \cite[\S 11.6]{Pi}). In our case,  $S$ is commutative.
The group  $G=A^*$  on the invertible elements of  $A$ is a finite group of triangular type    $G=H+J$, where  $H=S^*$. If  $A$ is the algebra of triangular matrices, then $G=\Bn=\TnF$ is the triangular group. The supercharacter theory for  $G=A^*$ was constructed in  \cite{P1}.

 Let $G=H+J$ be a finite group of triangular type.  Construct the group  $\tG$ that consists of the triples  $\tau=(t,a,b)$, where $t\in H$,~ $a,b\in N$, with operation   $$(t_1,a_1,b_1)\cdot (t_2,a_2,b_2) = (t_1t_2,~ t_2^{-1}a_1t_2a_2, ~ t_2^{-1}b_1t_2b_2).$$
The group   $\tG$ acts on   $J$  by the formula
\begin{equation}\label{rororo}
 \rho_\tau(x) = taxb^{-1}t^{-1}.
 \end{equation}
The representation of the group $\tG$ in the dual space $J^*$ is defined as usual  $$\rho^*_\tau\la(x) = \la(\rho(\tau^{-1})(x)).$$   In the space  $J^*$, there are defined also the left and right linear actions of the group  $G $ by the formulas  $b\la(x)= \la(xb)$ and $\la a(x) = \la(ax)$. Then $\rho_\tau(\la) =
tb\la  a^{-1}t^{-1}$.

 For any idempotent  $e\in kH$, we denote by  $A_e$ the subalgebra $eAe$.
 The subalgebra  $J_e=eJe\subset J$ is a radical in $A_e$. Denote  $e'=1-e$. The decomposition
 $$J= eJe\oplus eJe'\oplus e'Je\oplus e'Je'$$
 is called the Pierce decomposition.
 The dual space  $J_e^*$ is naturally identified with the subspace in  $J^*$ that consists of linear forms equal to zero on all components of Pierce decomposition except for the first component.

 Observe that, since the group   $H$ is abelian,  $he=eh=ehe$ for all  $h\in H$. The subset $H_e=eHe$ is a subgroup in the group of all invertible elements  of subalgebra  $A_e$. The subgroup $G_e=eGe=H_e+J_e$ is a finite group of triangular type, and this group is associated with the algebra  $A_e$ in the same way as   $G$ is associated with   $A$.
We define the group  $\tG_e$ similar to  the group  $\tG$. The map $h\to he$ is a homomorphism of the group  $H$ onto $H_e$ with the kernel
\begin{equation}\label{he} H(e) =\{h\in H: ~ he=e\}.\end{equation}

 The following definition is equivalent to the corresponding definition from  \cite{P1,P2}, although it differs by a form.\\
 \Def\Num\label{sing}. An $\rho_{\tG}$-orbit $\Oc$ is  \emph{singular} (with respect to  $H$) if $\Oc\cap J_e\ne\varnothing$ for some idempotent  $e\ne 1$ from $kH$. Otherwise,  the orbit $\Oc$  is \emph{regular} (with respect to  $H$).
   The elements of the singular (regular) orbits are called singular (regular) elements. Similarly, we define the singular and regular orbits  and elements in  $J^*$.

  The subgroup $H(e)$ admits the following characterization.\\
\Prop\Num\label{charac}~\cite[Lemma 2.5]{P2}.
1)~  $H(e)   = H_{y,\mathrm{right}}\cap H_{y,\mathrm{left}}$ for any regular (with respect to $H_e$) element  $y\in J_e$.\\
2)~ $H(e) = H_{\la,\mathrm{right}}\cap H_{\la,\mathrm{left}}$ for any regular (with respect to  $H_e$) element  $\la\in J_e^*$.

In the papers \cite{P1,P2}, the following statements are proved for any  $\tG$-orbit $\Oc$ in $J$. \\
 \Prop\Num\label{state}. 1) For any idempotent $e\in kH$, the intersection $\Oc\cap J_e$ is a $\tG_e$-orbit in $J_e$.\\
 2) There exists a unique idempotent $e\in kH$ such that $\Oc\cap J_e$ is a regular  $\tG_e$-orbit in  $J_e$ (with respect to  $H_e$). \\
Similar statements are true for   $\tG$-orbits in   $J^*$.

It is known that, for any representation of the finite group in the finite dimensional linear space   $V$ defined over a finite field, the number of orbits in  $V$ and $V^*$ are equal \cite[Lemma 4.1]{DI}. This statement and the  properties of  $\tG$-orbits imply  that the number of regular (singular)  $\tG$-orbits an  $J$ and $J^*$ are equal \cite[Proposition 2.11]{P1}.

Turn to definition of superclasses in the group  $G$. For any $g\in G$ and $(t,a,n)\in\tG$ consider the element
\begin{equation}\label{rtau}
R_\tau(g) = 1+ ta(g-1)b^{-1}t^{-1}
\end{equation}
from the algebra $A=kH+J$.
  If $g=h+x$, then
$ R_\tau(g)=h\bmod J$. Therefore,  $ R_\tau(g)\in G$.
The formula  (\ref{rtau}) defines an action of the group  $\tG$ on $G$. \\
\Def\Num. We refer to the $\tG$-orbits in  $G$ as  \emph{superclasses}.

 The group  $G$ decomposes into superclasses.
Denote by $\Bx$ the set of all triples  $\beta = (e,h, \omega)$, where
 $e$ is an idempotent from  $kH$, ~ $h\in H(e)$, and $\omega$ is a regular (with respect to $H_e$) ~ $\tG_e$-orbit in $J_e$. All elements  $h+\omega$ are contained in the common superclass \cite[corollary 3.2]{P2}; we denote it by  $K_\beta$.\\
\Theorem\Num~ \cite[теорема 3.3]{P2}. The correspondence  $\beta\to K_\beta$ is a bijection betweem the set triples  $\Bx$ and the set of superclasses in $G$.

Denote by  $\Ax$ the set of triples $\al = (e,\theta, \omega^*)$, where
 $e$ is an idempotent in   $kH$, ~ $\theta $ is a linear character (one dimensional representation) of the group  $ H(e)$,  and $\omega^*$ is a regular (with respect to $H_e$)~  $\tG_e$-orbit in  $J^*_e$.
 Since the subgroup  $H(e)$ is abelian, then the number of its linear characters equals to the number of elements. The number of regular (with respect to  $H_e$)  $\tG_e$-orbits in $J_e$ and $J_e^*$ are common. Hence, $|\Ax|= |\Bx|$.

We turn to construction of supercharacters.
Let $\al = (e,\theta, \omega^*)\in \Ax$, choose $\la\in \omega^*$.

 Consider the subgroup   $G_\al = H(e)\cdot N_{\la,\mathrm{right}}$, where $N_{\la,\mathrm{right}}$ is the stabilizer of  $\la$ for the right action of the group  $N=1+J$  on  $J^*$.
The subgroup  $N_{\la,\mathrm{right}}$ is an algebra subgroup; it is presented in the form $N_{\la,\mathrm{right}}= 1+J_{\la,\mathrm{right}}$, where $J_{\la,\mathrm{right}}$ is the right stabiliser of  $\la$ in $J$.
The group  $G_\al$ is a sum
 \begin{equation}\label{Gal}
   G_\al = H(e) + J_{\la,\rt};
 \end{equation}
it is finite group of triangular type.

Fix a nontrivial character  $t\to  \varepsilon^{t}$ of the additive group of the field
$\Fq$ with values in the multiplicative group  $\Cb^*$. By the triple  $\al = (e,\theta, \omega^*)$ and $\la \in\omega^*$, we define the linear character of the group  $G_\al$ by the formula
\begin{equation}\label{thela}
    \xi_{\theta,\la}(g) = \theta(h)\varepsilon^{\la(x)},
\end{equation}
where $g=h+x$,~ $h\in H(e)$ and $x\in J_{\la,\mathrm{right}}$.  Let us show that   $\xi=\xi_{\theta,\la}$ is really a linear character:
$$\xi(gg') = \xi((h+x)(h'+x')) = \xi(hh' + h'x + x'h  + xx') =$$
$$\theta(hh')\varepsilon^{\la(h'x)} \varepsilon^{\la(x'h)}\varepsilon^{\la(xx')} = \theta(h)\theta(h')\varepsilon^{\la(x)} \varepsilon^{\la(x')} = \xi(g)\xi(g').$$
The induced character
\begin{equation}\label{Indchi}
\chi_\al = \Ind(\xi_{\theta,\la}, G_\al, G)
\end{equation}
is called the  {\it supercharacter}.\\
\Theorem\Num ~\cite[Proposition 4.1, 4.4]{P1,P2}. \\
1) The supercharacters  $\{\chi_\al:~ \al \in\Ax\}$ are pairwise disjoint;\\
 2) Each supercharacter  $\chi_\al$ is constant on each superclass  $K_\beta$;\\
 3) $\{1\}$ is a superclass   $K(g)$ for  $g=1$.

 Applying the Proposition   \ref{dilemma},  we conclude.\\
 \Theorem\Num ~\cite[Theorem 4.5]{P2}. The systems of supercharacters $ \{\chi_\al|~ \al\in \Ax\}$ and superclasses   $\{\Kc_\beta|~ \beta\in \Bx\}$  give rise to  a supercharacter theory on the group   $G$.\\
 \Remark. While constructing the supercharacters, it is not  easy to find idempotents in the group algebra $kH$.  One can prove ~\cite[section 5]{P2} that, in the example \ref{exthree}, we may choose idempotents in the subalgebra  $S$ (in the case  $G=\mathrm{T}(n,\Fq)$, the subalgebra $S$ coincides with the subalgebra of diagonal matrices), and, in the example \ref{extwo}, we may  confine to  two idempotents   $ \left(\begin{array}{cc}0&0\\0&1\end{array}\right)$ and $ \left(\begin{array}{cc}1&0\\0&1\end{array}\right)$.

Consider the problems of restriction and induction in the constructed supercharacter theory \cite{P3}.
Let $G=H+J$ be a finite group of triangular type.
 Let $H$ be a subgroup in  $H$, ~ and  $J'$ be a subalgebra in  $J$ invariant with respect the left-right action of the group   $H'\times H'$  on $J$. Then $G'=H'+J'$ is a subgroup in  $G$; we call it the  \emph{subgroup of triangular type} in $G$.\\
\Theorem\Num\label{restr}~ \cite{P3}. The restriction of  supercharacter of the group  $G$ on its subgroup of triangular type is a sum supercharactes of this subgroup with nonnegative integer coefficients.

 Let $\phi $ be a superclass function (i.e., the complex valued function constant on superclasses)  on $G'$. Denote by  $\dot{\phi}$ the function on  $G$ equal to  $\phi$ on $G'$ and zero outside  $G'$.

Define the superinduction as follows:
 $$ \SInd\, \phi(g) = \frac{|H|}{|G|\cdot|G'|}\sum_{\tau\in \tG} \dot{\phi}(\rho(\tau)(g))=
\frac{|H|}{|G|\cdot|G'|}\sum_{a,b\in N,~ t\in H} \dot{\phi}(1+ta(g-1)bt^{-1}).$$
Easy to see that $\SInd\,\phi$ is a superclass function on  $G$. The supercharacter analog of the Frobenius theorem is valid. \\
\Theorem\Num\label{last}. Let $\psi$ be a superclass function on $G$. Then  $(\SInd\,\phi, \psi) = (\phi, \Res\,\psi)$.

    Let $\{\chi_\al\}$ be the system of supercharacters for the finite group of triangular type $G=H+J$,
   and $\{\phi_\eta\}$ be the system of supercharacters of its subgroup of triangular type $G'=H'+J'$. By theorem  \ref{restr}, $$\Res\,\chi_\al=\sum_\eta m_{\al,\eta}\phi_\eta, ~ \quad  m_{\al,\eta}\in\Zb_+ .$$
    The system of supercharacters form a basis in the space of superclass functions.
   Since  $\SInd\,\phi_\eta$ is a superclass function on the group  $G$, we have $$ \SInd\,\phi_\eta=\sum_{\al}
 a_{\eta,\al}\chi_\al. $$
 The theorem \ref{last} implies $$ a_{\eta,\al} = \frac{m_{\al,\eta}(\phi_\eta,\phi_\eta)}{( \chi_\al, \chi_\al)}.$$
 \Cor\Num. For any supercharacter  $\phi$ of the subgroup of triangular type $G'$,  the superinduction  $\SInd\,\phi$ is a sum of supercharacters of the group  $G$ with nonnegative rational coefficients.

 Our next goal is to present the supercharacter version of the A.A.Kirillov formula in the case of the finite groups of triangular type. This formula generates the formulas  (\ref{kirdi}) and (\ref{kirdisec}) for the algebra groups.

Any superclass  $K$ finite group of triangular type $G$  has the element  of the form $g=h+x$, where   $hx=xh=x$. The element   $h-1\in kH$ is associated with some idempotent  $f'\in kH$ (i.e., these elements differ by an invertible  multiplier  from $kH$).  Let $f=1-f'$. The condition $hx=xh=x$ is equivalent to  $x\in J_f$.  Two elements $g=h+x$ and $g'=h+x'$, where $ x,~ x'\in J_f$, belong to a common superclass whenever  $x$ and  $x'$ belong to a common $\tG_f$-orbit \cite[Theorem 3.1]{P2}.

 Since supercharacters are constant on superclasses, it is sufficient calculate the values of supercharacters on the elements of the form   $g=h+x$, where   $hx=xh=x$.
  Let $\al= (e,\theta,\omega^*)$, and    $\chi_\al$  be a supercharacter.
If  $h\notin H(e)$, then  $g$ does not belong to the  $\Ad_G$-orbit of the subgroup  $G_\al$. By definition of the induced character, we obtain $\chi_\al(h)=0$.

  Consider the case $h\in H(e)$, i.e.,  $he=e$; then  $f'e=0$ and  $e<f$.
We see $J_e\subset J_f $, and    $J_e^*$ is a subspace in  $J_f^*$.
In particular,  $\omega^*\in J_f^*$. There exists a unique $\rho^*(\tG_f)$-orbit  $\Omega^*$ in $J_f^*$ such that its intersection with  $J_e^*$ coincides with  $\omega^*$ (see Proposition  \ref{state}).
Denote by   $\dot{\theta}(h)$  the function   $H\to\Cb$ equal to $\theta(h)$ on $H(e)$ and  zero ontside  $H(e)$.\\
   \Theorem \Num\label{kirillov}.
The value of the supercharacte $\chi_\al$ on the element  $g=h+x$, where $hx=xh=x$, can be calculated by the formula
  \begin{equation}\label{kirform1}
\chi_\al(g) = \frac{|H_e|\cdot\dot{\theta}(h)}{n(\Omega^*)} \sum_{\mu\in \Omega^*} \varepsilon^{\mu(x)},
\end{equation}
where $n(\Omega^*)$ is the number of right  $N_f$-orbits in $\Omega^*$.\\
The other form of A.A.Kirillov formula has the form.  \\
\Theorem\Num. The value of the supercharacter $\chi_\al$ on the element  $g=h+x$, where $hx=xh=x$, can be calculated by the formula
\begin{equation}\label{kirform}
\chi_\al(g) = \frac{|H_e|\cdot |\la N_f|\cdot\dot{\theta}(h)}{|\rho(\tG_f)(x)|} \sum_{y\in \rho(\tG_f)(x) } \varepsilon^{\la(y)}.
\end{equation}

 The next goal is  to obtain the group of triangular type version of the theorem \ref{intersec}.   The scalar product on the group   $G$ is defined as in  (\ref{scalar}).\\
 \Lemma\Num\label{jjj}. Let $\la\in J_f^*$. Then  $J_f\la\cap \la J_f = J\la\cap \la J$.\\
 \Proof. Since $\la\in J_f^*$, we have $\la f=f\la=\la$. The  inclusion  $J_f\la\cap \la J_f\subseteq J\la\cap \la J$ is obvious. Let  $\mu\in J\la\cap \la J$.  Then $\mu=j_1\la$ and $\mu=\la j_2$ for some  $j_1, j_2\in J$.
  Hence,  $\mu f = j_1\la f= j_1\la =\mu$ and  $f \mu = f\la j_2 = \la j_2 = \mu$.
  Therefore,  $$\mu = f\mu f\in f(J\la\cap \la J)f=(fJf)(f\la f)\cap (f\la f)(fJf) = J_f\la\cap \la J_f. $$
  $\Box$\\
 \Theorem\Num. Let $\al=(e,\theta, \omega^*)\in \Ax$ and $\la\in \omega^*$. Then
  \begin{equation}\label{scachial}
(\chi_\al,\chi_\al) = \frac{|H_{N\la N}|}{|H(e)|}\cdot |J\la\cap \la J|,
\end{equation}
where $H_{N\la N}$ is the stabilizer of $N\la N$ with respect to  $\Ad_H^*$ action the subgroup  $H$.\\
\Proof. The subgroup $G'=H(e) N$ contains  $G_\al$. Consider the character
$\chi'_{\theta,\la}=\Ind(\xi_{\theta,\la},G_\al, G')$ of the subgroup $G'$.  Applying the Intertwining Number Theorem  \cite[теорема 44.5]{CR}, we obtain
$$(\chi_\al,\chi_\al)=\left(\Ind (\chi'_{\theta,\la},G',G), \Ind (\chi'_{\theta,\la},G',G)\right) =
\sum_{h\in H/H(e)} (\chi'_{\theta,\Ad^*_h\la}, \chi'_{\theta,\la}),$$
The subgroup  $G'$ is also a finite subgroup of triangular type; its diagonal part   $H(e)$ identically acts on  $\la$. Two characters  $ \chi'_{\theta,\Ad^*_h\la}$ и  $\chi'_{\theta,\la}$ are disjoint if  $\Ad^*_h\la\notin N\la N$.
Therefore,
\begin{equation}\label{firstchi}
(\chi_\al,\chi_\al)=\frac{|H_{N\la N}|}{|H(e)|}(\chi'_{\theta,\la},\chi'_{\theta,\la}).
\end{equation}

We apply the formula  (\ref{kirform1})  to the supercharacter  $\chi'_{\theta,\la}$ of the group  $G'$.   For
$g=h+y\in G'$,~ $hy=yh=y$, we have
$$\chi'_{\theta,\la}(g) = \theta(h) \chi_{\theta,\la}^f(1+y),$$
 where $\chi_{\theta,\la}^f$ is a supercharacter of the group  $N_f$ constructed by  $\la\in J_e^*\subseteq J_f^*$.

Fix  $h\in H(e)$ and  $x\in J$ such that  $hx=xh=x$.
As  above,  the condition on  $x$ is equivalent to $x\in J_f$,  where $f=1-f'$ and the idempotent  $f'$ is associated to  $h-1$.
 Each element of the superclass  $K(h+x)$  in  $G'$ can be uniquely presented in the form   $(1+v)(h+y)(1+u)$, where $v\in fJf'$, $u\in f'J$, and $y\in N_fxN_f$.
The supercharacter is constant on superclasses; then
$$
  \sum_{g\in K(h+x)}\chi'_{\theta,\la}(g)\overline{\chi'_{\theta,\la}(g)}= \frac{|J|}{|J_f|}\sum_{y\in N_fxN_f}  \chi_{\theta,\la}^f(1+y)\overline{ \chi_{\theta,\la}^f(1+y)} .
$$
Hence
  \begin{multline}\label{sumsum}
    \sum_{g\bmod J=h }\chi'_{\theta,\la}(g)\overline{\chi'_{\theta,\la}(g)} = \frac{|J|}{|J_f|}  \sum_{y\in J_f} \chi_{\theta,\la}^f(1+y)\overline{ \chi_{\theta,\la}^f(1+y)} = \\ |J| \cdot (\chi_{\theta,\la}^f,\chi_{\theta,\la}^f) = |J| \cdot |J_f\la\cap\la J_f| = |J| \cdot |J\la\cap\la J|.
     \end{multline}
Observe that the sum  (\ref{sumsum}) does not depent on  $h\in H(e)$. Then
\begin{multline}\label{secondchi}
(\chi_\al,\chi_\al)=\frac{|H_{N\la N}|}{|H(e)|}\cdot\frac{1}{|G'|}\sum_{g\in G'}\chi'_\al(g)\overline{\chi'_\al(g)}=\\
\frac{|H_{N\la N}|}{|H(e)|\cdot|G'|}\cdot |H(e)|\cdot|J| \cdot |J\la\cap\la J| = \frac{|H_{N\la N}|}{|H(e)|} \cdot |J\la\cap\la J|.\quad \Box
\end{multline}

\subsection{The supercharacter theory for the triangular group}

Consider the algebra  $A=\tx(n,\Fq)$ that consists of all   $n\times n$-matrices
 with entries  from the field  $\Fq$ and zeros below the diagonal. The {\it triangular group}  $G=\sB_n=\TnF$ is a group of all invertible elements in  $A$; it is a finite group of triangular type  $G=H+J$, where   $J$ is a subalgebra $\ut_n=\utn$ of all triangular matrices with zeros on the diagonal, and $H=\{(a_1,\ldots, a_n):~ a_i\in\Fq^*\}$ is the subgroup of diagonal matrices. The subgroup $G$ is a semidirect product  $G=HN$, where $N=\UT_n$. In this subsection, we concretize the constructed supercharacter theory to the triangular group.

As above, a positive root is a  pair  $(i,j)$, where $1\leq i<j\leq n$.
 Let $D$ be a basic subset (see  subsection \ref{basic}). By a \emph{support} of  $D$, we call  the subset $\supp (D)=\row(D)\cup\col(D)$ in $[n]$. For each basic subset $D$, we construct  the element $$x_D=\sum_{(i,j)\in D}  E_{ij}$$ in $J$.
Denote by  $\Oc_D$ the orbit of the element  $x_D$ with respect to the action $\rho$ of the group  $\tG$ on $J$ (see (\ref{rororo})). Each $\tG$-orbit in  $J$ has the form  $\Oc_D$ for some basic subset  $D$. Similarly for  $J^*$; each $\tG$-orbit in $J^*$ is the orbit    $\Oc_D^*$  of some $$
 \la_D=\sum_{(i,j)\in D}  E^*_{ij}.$$
\Lemma\Num\label{nnn}~\cite{P1}. The orbit $\Oc_D$ (respectively, $\Oc^*_D$) is regular if and only if  $\supp(D)=[1,n]$.

By basic subset   $D$, we  construct idempotent
$$e=e_D=\sum_{i\in \supp(D)}E_{ii}.$$
The subgroup   $H$ decomposes into the product of subgroups  $H_i=\{(1,\ldots,a_i,\ldots 1)\}$ isomorphic to $\Fq^*$. The subgroups $H(e)$ and $H_e$, defined in (\ref{he}), have the form
$$ H_{e}=\prod_{i\in\supp(D)} H_i,\quad H(e)=\prod_{i\in [n]\setminus\supp(D)} H_i.$$
The subgroup $H$ is the product  $H=H_e\cdot H(e)$.

For the idempotent $e=e_D$, the subgroup $G_{e}$ is naturally isomorphic to triangular subgroup $\sB_m$, where $m$  is the number of elements in $\supp(D)$.
The element  $x_D$ belongs to  $J_{e}=eJe$ (respectively,  $\la_D$ belongs to $J_{e}^*$).
  Denote by  $\omega_D$  the $\tG_{e}$-orbit of
$x_D$  в $J_{e}$. Respectively,   $\omega^*_D$ is the $\tG_{e}$-orbit of $\la_D$ in $J_{e}^*$. By Lemma \ref{nnn}, we obtain that  $\omega_D$
and $\omega^*_D$ are regular orbits in $J_{e}$ and $J_{e}^*$ (with respect to $H_{e}$). \\
{\bf Example}. $G=\mathrm{T}(3,\Fq)=\left(\begin{array}{ccc} *&*&*\\
0&*&*\\0&0&*\end{array}\right)$, and  $D=\{(1,3)\}$. Then
 $$e=e_D=\left(\begin{array}{ccc} 1&0&0\\
0&0&0\\0&0&1\end{array}\right),~ H_{e}= \left(\begin{array}{ccc} *&0&0\\
0&1&0\\0&0&*\end{array}\right),~
 H(e)= \left(\begin{array}{ccc} 1&0&0\\
0&*&0\\0&0&1\end{array}\right),$$
$$ J_{e}= \left(\begin{array}{ccc} 0&0&*\\
0&0&0\\0&0&0\end{array}\right),\quad
 G_{e}= \left(\begin{array}{ccc} *&0&*\\
0&1&0\\0&0&*\end{array}\right),$$ $$ \omega_D= \left(\begin{array}{ccc} 0&0&\ne 0\\
0&0&0\\0&0&0\end{array}\right) \mathrm{is~~ a~~ regular}~~\tG_e-\mathrm{orbit~~in~~}  J_{e}.$$

Given  the basic subset  $D$ and the element  $h\in H(e)$, we construct the superclass  $K_{h,D'}$ as the $\tG$-orbit of the element  $h+x_{D}$ in the group  $G$. For $D$ and the character $\theta$ of the group  $H(e)$, we define the supercharacter  $\{\chi_{\theta,D}\}$ of the group $G$ as  in the subsection \ref{supertrtr} for the triple $(e,\theta,\omega)$.\\
\Theorem\Num\label{triang}~\cite{P1}. The systems of supercharacters  $\{\chi_{\theta,D}\}$ and superclasses   $\{K_{h,D}\}$ give rise to a supercharacter  theory  of the group  $G=\Bn$.

Let us calculate the values of the supercharactor   $\chi_{\theta,D}$ on the superclass
$K_{h,D'}$ of the element $h+x_{D'}$. We need some new notations.

For each positive root $\gamma=(i,j)$, ~$1\leq i<j\leq n$, we denote
$$ \Delta'(\gamma) =\{(i,k)| ~i<k<j\},\quad\quad\Delta''(\gamma) = \{(k,j)|~
i< k<j\}.$$
 The numbers of elements  of  both subsets a equal and coincide with $j-i-1$.
 We take
  $ \delta' (D,D') = 0$ if there exists $\gamma\in D$ and $\gamma\in D'$ such that  $\gamma'\in\Delta'(\gamma)$. Otherwise, we take  $\delta'(D,D')=1$.
Similarly, we define  $ \delta''(D,D')$.

 Take $\delta_0(D,h) = 1$ if  $h\in H(e_D)$. Otherwise,   $\delta_0(D,h) = 0$. Denote
\begin{equation}\label{delta}
    \delta(D,h,D') =  \delta'(D,D') \delta''(D,D')\delta_0(D,h).
\end{equation}

For any positive root  $\gamma=(i,j)$, ~$1\leq i<j\leq n$, we denote by
$P(\gamma)$ the submatrix of the matrix  $h-1+x_{D'}$ with systems of rows and columns
$[i+1,j-1]$.  Let  $m(\gamma, h, D')$ denote corank of  $P(\gamma)$.  Since there is at most one element in each row and column of  $P_\gamma$,  ~  $m(\gamma, h, D')$  is the number of zero rows (columns) in  $P_\gamma$.  Introduce the notations
\begin{equation}\label{corank}
    m(D, h, D') = \sum_{\gamma\in D} m(\gamma, h, D'),\quad \quad s(D,D') = |D|+|D\setminus D'|.
\end{equation}
\Theorem\Num\label{vvv}. The value of the supercharacter $\chi_{\theta,D}$ on the superclass
$K_{h,D'}$ equals to
\begin{equation}\label{value}
\chi_{\theta,D}(K_{h,D'})  = \delta(D,h,D')(-1)^{|D\cap D'|}q^{m(D, h, D')}(q-1)^{s(D,D')}\theta(h).
\end{equation}

\section{The Hopf algebra of supercharacters of the triangular group}\label{superhopftr}

\subsection{The Hopf algebra $\NPS$}  Let we have  the system of non-commuting variables $X\cup Y$
that  consists of two parts $X=\{x_1,x_2,\ldots\}$ and $Y=\{y_1,\ldots,y_{|Y|}\}$.  As in  subsection  \ref{hau}, the group $S_\infty=\varinjlim S_n$ acts on the set $X$ by finite permutations.
Extend this action to the set $X\cup Y$ setting  it  identical on $Y$.
  Consider the linear space  $\NPS$ of all
 formal complex power series in variables of $X\cup Y$ of   bounded degree in $X$ and finite degree in $Y$, and invariant with respect to the group  $S_\infty$.
The linear space
    $\NPS$ is a direct sum of the subspaces  $$ \NPS = \sum_{n=0}^\infty \NPSn $$
    of homogeneous  elements of fixed degree.  For $n=0$, the subspace  $\mathrm{NPS}_0(X,Y)$  coincides with the field $\Cb$.

The  linear space  $\NPS$ is a graded algebra under the natural multiplication. We call the algebra  $\NPS$ the \emph{algebra of partially symmetric functions in non-commuting variables}.

Let $A$ be  a  finite subset of positive integers.  \\
\Def\Num. By a \emph{rigged partition} of  $A$ we call a pair  $\Pc=(P,\Phi_\Pc)$, where $P$ is a partition of some subset  $\supp(P)\subseteq A$ (the partial partition of the set $A$), and $\Phi_\Pc$ is a map  $A\setminus \supp(P)\to [|Y|]$. Denote $\Pc\vDash A$.\\
\Def\Num.  Let  $\Pc\vDash A$ and $\Pc=(P,\Phi_\Pc)$. Then  $\Pc=\Pc_1+\Pc_2$, where $\Pc_i\vDash A_i$,~  $\Pc_i=(P_i,\Phi_{\Pc_i})$, if  $A=A_1\sqcup A_2$, ~ $P=P_1+P_2$ (see subsection \ref{hau}), the map  $\Phi_\Pc$ coincides with $\Phi_{\Pc_i}$  being reduced on $A_i\setminus \supp(\Pc_i)$.

Let  $\Pc=(P,\Phi_\Pc)$ be a ridded partition of the set $ [n]$.  Enumerate  the components of the partial partition  $P=\{\Pi_1,\ldots, \Pi_s\}$.
 Correspond $\Pc$ to the element  $m_\Pc\in\NPSn$ that is a sum of all monomials of the form  $\sigma (z_1\cdots z_n)$, where $\sigma\in S_\infty$, ~  $z_i=y_{\Phi(i)}$ for all  $i\in [n]\setminus \supp(P)$, and
$z_i=x_j$ for $i\in \Pi_j$. The element  $m_\Pc$ does not depend on the enumeration of  partial partition  $P$.\\
{\bf Example}. $n=4$, ~ $P=13|4$ и $\Phi(2)=1$. Then $m_\Pc$ is the sum $x_1y_1x_1x_2 + x_2y_1x_2x_1 + x_1y_1x_1x_3 +  \ldots$.

The system of elements $\{m_\Pc:~ \Pc\vDash [n]\}$ is a basis in $\NPSn$.
 For $\Pc\vDash [k]$ and $\Qc\vDash [n-k]$,  we consider the  rigged partition  $\Sc=(S,\Phi_\Sc)$, where $S=P\sqcup \{k+Q\}$, and $\Phi_\Sc(i)=\Phi_\Pc(i)$  if $1\leq i\leq k$, and $\Phi_\Sc(k+j)=\Phi_\Pc(j)$ if  $1\leq j\leq n-k$. Denote   $\Sc=(\Pc|\Qc)$.

We say that the rigged partition  $\Rc=(R,\Phi_\Rc)$ is a \emph{direct consequence} of the partition  $(\Pc|\Qc)$ if  $\supp(R)=\supp(P)\sqcup\supp(Q)$, $R\wedge ([k]|[n-k])=(P|Q)$, and $\Phi_\Rc$ coincides with  $\Phi_{(\Pc|\Qc)}$. Denote $ (\Pc|\Qc)\to\Rc$.
 Then
 \begin{equation}\label{mm}
   m_\Pc m_\Qc=\sum_{(\Pc|\Qc)\to\Rc} m_\Rc.
 \end{equation}
 We define a comultiplication  by the formula
 \begin{equation}\label{dede}
   \Delta(m_\Pc)=\sum_{\Pc=\Pc_1+\Pc_2}m_{\st(\Pc_1)}\otimes m_{\st(\Pc_2)}.\end{equation}
 The counit is defined as  the map that corresponds the  series to its constant term.\\
\Theorem\Num.  The algebra  $\NPS$ is a graded Hopf algebra under the defined comultiplication and counit. \\
\Proof. Since the space  $\NPS$ is a connected graded algebra (i.e.,  $\mathrm{NPS}_0(X,Y)= \Cb$), it is sufficient to prove that $\NPS$ is a bialgebra  (\cite{Zel,GR}), that is, the comultiplication is an algebra homomorphism. For  $\Pc\vDash [k]$ and $\Qc\vDash [n-k]$, applying  (\ref{mm}) and (\ref{dede}), we obtain

\begin{multline*}\label{fff}
   \Delta(m_\Pc) \Delta(m_\Qc) = \left( \sum_{\Pc=\Pc_1+\Pc_2}m_{\st(\Pc_1)}\otimes m_{\st(\Pc_2)}\right) \left( \sum_{\Qc=\Qc_1+\Qc_2}m_{\st(\Qc_1)}\otimes m_{\st(\Qc_2)}\right) = \\
 \sum_{\Pc=\Pc_1+\Pc_2,~~ \Qc=\Qc_1+\Qc_2}m_{\st(\Pc_1)}m_{\st(\Qc_1)}\otimes m_{\st(\Pc_2)} m_{\st(\Qc_2)} =
   \sum_{\Rc=\Rc_1+\Rc_2,~ (\Pc|\Qc)\to\Rc}m_{\st(\Rc_1)}\otimes m_{\st(\Rc_2)}.
\end{multline*}

On the other hand,
\begin{multline*}
 \Delta(m_\Pc m_\Qc) = \Delta\left(\sum_{(\Pc|\Qc)\to\Rc} m_\Rc\right) =  \sum_{\Rc=\Rc_1+\Rc_2,~~ (\Pc|\Qc)\to \Rc}m_{\st(\Rc_1)}\otimes m_{\st(\Rc_2)}.
  \end{multline*}
  Hence $ \Delta(m_\Pc m_\Qc) = \Delta(m_\Pc) \Delta(m_\Qc)$.~ $\Box$

\subsection{The Hopf algebra   $\SCB$}

Let   $\Bn = \UTn$. Consider the space of superclass functions  $\SCB_n$ on the group $\Bn$ and the direct sum  $$ \SCB=\sum_{n=0}^\infty \SCB_n.$$
The multiplication and comultiplication in $\SCB$ are defined below similarly to the algebra $\SCU$.

The group $\sB_k\times \sB_m$ is a factor group of  $\sB_n$, where $n=k+m$.
Let  $\phi\in \SCB_k$ and $\psi\in \SCB_m$, then, by definition,     \begin{equation}\label{muB}
\phi\cdot \psi = \Inf(\phi\times\psi)\in \SCB_n,
\end{equation}
   where the inflation is a composition of  $\phi\times\psi$ and the natural projection $\sB_n\to \sB_k\times \sB_m$.

 For any partition  $T=(A_1|A_2)$ of the set $[n]$ , there exists the subgroup $\sB_{A_1}\times \sB_{A_2}$ of the group  $\sB_n$; the subgroup consists of all unitriangular matrices such that  $a_{ij}\ne 0$ implies that  $i$, $j$  lie in a common component of the partition $T$.  The subgroup  $\sB_{A_1}\times \sB_{A_2}$ is naturally isomorphic to $\sB_{|A_1|}\times \sB_{|A_2|}$ (denote the isomorphism by  $\pi_T$) .
 Let $\chi\in \SCB_n$. Then, by definition,
  \begin{equation}\label{coB}
 \Delta(\chi) = \sum_{T=(A_1|A_2)\vdash [n]} { ^T\Res}^{\sB_n}_{\sB_{|A_1|}\times \sB_{|A_2|}}(\chi),\end{equation}
 where $${ ^T\Res}^{\sB_n}_{\sB_{|A_1|}\times \sB_{|A_2|}}(\chi)(g) =\Res^{\sB_n}_{\sB_{|A_1|}\times \sB_{|A_2|}} (\chi)(\pi_T^{-1}(g)).$$

 Suppose that $|Y|=q-2$. Fix enumerations in the set $\Fq^*\setminus \{1\}=\{u_1,\ldots,u_{q-2}\}$ and in  set  $\Irr(\Fq^*)\setminus \{1\}=\{\eta_1,\ldots,\eta_{q-2}\}$.
 We correspond the rigged partition  $\Pc=(P,\Phi_\Pc)\vDash [n]$  to  the supercharacter as follows.
 By the partial partition $P$, we construct the basic subset  $D_P$ in the set of positive roots, as in subsection  \ref{scu}. Observe that $\supp(P)\subseteq \supp(D)$ and the equality holds if the partition  $P$ has no singleton components.

 Construct the element  $$h_\Pc=\prod_{i\in [n]\setminus \supp(P)} u_{\Phi_\Pc(i)}\in H(e_{D_P})$$ and the superclass $K_\Pc$  of the element  $h_\Pc+x_{D_P}$.   We constrct also the supercharacter  $\chi_{\Pc}$ by the subset  $D_P$ and by the character  $$\theta_\Pc=\prod_{i\in [n]\setminus \supp(P)} \eta_{\Phi_\Pc(i)}$$ of the subgroup   $H(e_{D_P})$.
 In terms of rigged partitions, Theorem  \ref{triang}  can be presented in the following way.\\
 \Prop\Num. The systems of supercharacters  $\{\chi_\Pc\}$ and superclasses   $K_{\Pc}$, where $\Pc$ runs through the set of all rigged partitions of  $[n]$, give rise  to a supercharacter theory of the triangular group  $\sB_n$.

 We denote by $\kappa_\Pc$ the characteristic function of the superclass  $K_\Pc$. Let as above $|Y|=q-2$.\\
\Theorem\Num.  The linear space  $\SCB$ is a Hopf algebra  with respect to the defined multiplication and comultiplication. The map  $\kappa_\Pc\to m_\Pc$, where   $\Pc$ runs through the set of all rigged partitions of  $ [n]$,  can be linearly extended to the isomorphism of the Hopf algebra  $\SCB$ onto the Hopf algebra  $\NPS$.\\
\Proof.   It is sufficient to prove the  multiplication and comultiplication of the elements $\kappa_\Pc$ satisfies the equalities  (\ref{mm}) and (\ref{dede}).

The multiplication is defined by the formula (\ref{muB}). Let  $e_{[k]}$ be the idempotent equal to the sum of first  $k$  diagonal matrix units, and  $e_{[n-k]}$ be the sum of the last $n-k$ ones. Then $\kappa_\Pc\kappa_\Qc$ is the characteristic function of the set  $$\Lambda = \mathrm{diag}(K_\Pc|K_\Qc)+ e_{[k]}Je_{[n-k]}.$$

 Denote $G=\sB_n$. The set $\Lambda$ is invariant with respect to the action of the group  $\tG$ and, therefore, it decomposes into superclasses (see subsection  \ref{supertrtr}). All elements from  $\Lambda$  has common diagonal part
$h=\mathrm{diag}(h_\Pc,h_\Qc)$. Let $f'$ be  the idempotent associated with the diagonal matrix $h-1$. For the idempotent  $f=1-f'$, there defined the subgroup  $G_f=H_f+J_f$, where $H_f=fHf$ and $J_f=fJf$. The subgroup  $G_f $ is isomorphic to  $\sB_s$
for some  $s\leq n$; it is a semidirect product  $G_f=H_fN_f$, where $N_f=1+J_f$ is isomorphic to  $\UT_s$. Similarly for  $h_\Pc$ and $h_\Qc$ one can construct the idempotents $f_\Pc\leq e_{[k]}$ and $f_\Qc\leq e_{[n-k]}$. Then $f=f_\Pc+f_\Qc$;  $N_f$ contains the subgroup
$N_{f_\Pc}\times N_{f_\Qc}$ isomorphic to $\UT_{s_1}\times \UT_{s_2}$, where $s=s_1+s_2$.

Each superclass of  element with the diagonal part  $h$ contains the subset  $h+\omega$, where $\omega$ is some orbit in  $J_f$ with respect to the group  $\tG_f$ (see subsection  \ref{supertrtr}).
So,   the set $1+\omega$  is the   $H_f$-orbit in the set of superclasses of the group    $N_f=\UT_s$; it is uniquely determined by the given superclass with diagonal part  $h$.
If, in addition, this superclass belongs to  $\Lambda$, then the projection of
 $1+\omega$ on  $\UT_{s_1}\times \UT_{s_2}$
   coincides with the $H_f$-orbit of the superclass of the element $\mathrm{diag}(1+x_{D_P},1+x_{D_Q})$.
  Applying statements of subsection  \ref{huni}, $1+\omega$ corresponds to the partition $R$ of the set $[s]$ and  $R\wedge ([s_1]|[s_2])=(P|Q)$ (see Remark  \ref{hopfave}). The superclasses that belong to   $\Lambda$ are the superclasses of the group $\sB_n$ of the form $K_\Rc$, where  $(\Pc|\Qc)\to \Rc$. This proves
 \begin{equation}
   \kappa_\Pc  \kappa_\Qc=\sum_{(\Pc|\Qc)\to\Rc}  \kappa_\Rc.
 \end{equation}

According to the formula (\ref{dede}),  $ \Delta( \kappa_\Pc)$ is the sum of the restriction of $\kappa_\Pc$ on the subgroups $B_{A_1}\otimes B_{A_2}$.
 Each restriction is  the characteristic function of the intersection of the superclass $K_\Pc$ with the subgroup  $B_{A_1}\otimes B_{A_2}$ isomorphic to  $B_{|A_1|}\times B_{|A_2|}$.
  This intersection decomposes into superclasses $\{K_{\Pc_1}\times K_{\Pc_2}\}$ of the group  $B_{A_1}\otimes B_{A_2}$.
The $\tG$-orbit of each  superclass $K_{\Pc_1}\times K_{\Pc_2}$  coincides with  $K_\Pc$.
This proves that the intersection consists of a single superclass  and $\Pc=\Pc_1+\Pc_2$.
Hence
 $$
   \Delta( \kappa_\Pc)=\sum_{\Pc=\Pc_1+\Pc_2}\kappa_{\st(\Pc_1)}\otimes \kappa_{\st(\Pc_2)}.\quad \Box$$

\subsection{The dual Hopf algebra $\SCB^*$}

The linear space   $\SCB^*=\sum_{n=0}^\infty \SCB_n^*$ is a Hopf algebra with respect to operations dual to operations in  $\SCB$. By the scalar product, we identify   $\SCB^*$ and $\SCB$ as a linear spaces.
The definition of superinduction for the groups of triangular type enables to calculate the dual operations in term of the space $\SCB$ (see Theorem \ref{last}).

 Let $\phi\in \SCB_k$, ~ $\psi\in \SCB_m$ and $n=k+m$. For the partition  $T=(A|A^c)$ the dual operation to $^T\Res$ is the operation  $^T\SInd: \SCB_k\otimes \SCB_m\to \SCB_n$ defined as
$$^T\SInd(\phi\times\psi)= \SInd (\phi\pi_A\times\psi\pi_A).$$
 The multiplication in $\SCB^*$ is defined by the formula
 \begin{equation}\label{bsind}
\phi\cdot\psi = \sum_{T=(A_1|A_2),~ |A_1|=k, ~|A_2|=n-k} {^T\SInd}^{\mathrm{B}_n}_{\mathrm{B}_{|A_1|}\times \mathrm{B}_{|A_2|}}(\phi\times\psi).
\end{equation}
By definition,  comultiplication in  $\SCB^*$ is the deflation  $\Defl$ that is a  dual operation to the inflation $\Inf$. Let $\tau$ be the natural projection  $\mathrm{B}_n\to \mathrm{B}_k\times\mathrm{B}_m$ and  $\chi\in \SCB_n$. Then
\begin{equation}\label{bdefl}
\Defl(\chi) (u,v) = \frac{1}{| \tau^{-1}(1)|}\sum_{x\in\tau^{-1}(u,v)}\chi(x),
\end{equation}
where $u\in \mathrm{B}_k$, ~ $v\in\mathrm{B}_m$.\\
\Theorem\Num. The operations  (\ref{bsind}) and  (\ref{bdefl}) define  the structure of dual Hopf algebra $\SCB^*$ on the linear space $\SCB$.

The system  $\{\kappa_\Pc:~ \Pc\vDash [n]\}$ is a basis in  $\SCB_n$.
The dual basis  $\{\kappa_\Pc^*:~ \Pc\vDash [n]\}$ consists of the elements    $\kappa_\Pc^*=z_\Pc\kappa_\Pc$, where
$z_\Pc=|\mathrm{B}_n|/|K_\Pc|.$
Calculate the operations of  $ \SCB^*$ on the basic elements.
For $\Pc\vDash[k]$ and $\Qc\vDash[m]$ we have
$$\kappa_\Pc^*\cdot\kappa_\Qc^* = \sum \kappa^*_{(\st_{A_1}^{-1}(\Pc)|\st_{A_2}^{-1}(\Qc))},$$
where the sum is taken over all partitions  $(A_1,A_2)\vdash [n]$ such that  $|A_1|=k$, ~$|A_2|=n-k$.

 Let  $\Pc\vDash[n]$.  Then
  $$\Delta(\kappa_\Pc^*) = \sum_{k=0}^n \kappa^*_{\Pc_{[k]}}\otimes \kappa^*_{\Pc_{[k]^c}},$$
 where $\Pc_{[k]}$ and $\Pc_{[k]^c}$ are the intersections  of the rigged partition  $\Pc$ with  $[k]$ and its complement  $[k]^c$ in  $[n]$.


\begin{thebibliography}{99}

\bibitem{VERY}
 Aguiar M.,  Andr\`{e} C.,   Benedetti C., Bergeron N., Zhi Chen,  Diaconis P.,  Hendrickson A.,  Hsiao S.,  Isaacs I.M.,  Jedwab A.,  Johnson K., Karaali G.,  Lauve A., Tung Le,  Lewis S., Huilan Li,  Magaarg K.,  Marberg E.,   Novelli J-Ch., Amy Pang,  Saliola F.,  Tevlin L., Thibon  J-Y.,  Thiem N.,  Venkateswaran V., Vinroot C.R., Ning Yan,  Zabricki M.,  \emph{Supercharacters, symmetric functions in noncommuting variables, and related Hopf algebras}, Advances in Mathematics, 2012, vol. 229, no.4,  2310-2337.

\bibitem{A1}
Andr\'{e}~C.A.M., \emph{Basic characters of the unitriangular group}, J. Algebra, 1995, vol. 175, 287-319.

 \bibitem{A2}
Andr\'{e}~C.A.M., \emph{Basic sums of coadjoint orbits of the unitriangular group}, 1995, J. Algebra,
vol. 176,  959-1000.

 \bibitem{A3}
Andr\'{e}~C.A.M., \emph{The basic character table of the unitriangular group},  J. Algebra, 2001,
 vol.  241, 437-471.

 \bibitem{A4}
Andr\'{e}~C.A.M., \emph{The basic characters  of the unitriangular group(for arbitrary primes}, Proc. Am. Math. Soc.,  2002, vol. 130, 1943-1954.
 \bibitem{A5}
 Andr\'{e}~C.A.M.,\emph{ Basic sums of Coadjoint Orbits of the Unitriangular group}, J. Algebra, 1995,
 vol. 176, 959-1000.

 \bibitem{A6}
Andr\'{e}~C.A.M., \emph{Hecke algebras for the basic characters of the unitriangular group}, Proc.Amer.Math.Soc., 2003, vol.132, 987-996.

 \bibitem{A7}
Andr\'{e}~C.A.M., Nicol\`{a}s A., \emph{Supercharacters of the adjoint group of a finite radical ring}, J.Group Theory, 2008, vol. 11, 709-746.

 \bibitem{Asite}
 Andr\'{e}~C.A.M., \emph{ Supercharacters of unitriangular groups and set-partition combinatorics}, CIMPA school: Modern Methods in Combinatorics ECOS 2013,  San Luis, Argentina.

\bibitem{AN}
Andr\'{e}~C.A.M., Neto~A.M., \emph{A supercharacter theory for the Sylow p-subgroups of the finite symplectic and orthogonal groups},
J. Algebra, vol. 322, no. 4, 2009, 1273-1294.

\bibitem{Walk}
Arias-Castro~E., Diaconis~P., Stanley~R., \emph{A super-class walk on upper-triangular matrices}, J. Algebra, 2004,
vol. 278, 739-765.

\bibitem{Ash}
Asharafi A.R., Koorepazan-Moftakhar F., \emph{Towards the classification of finite simple groups with exactly three or four supercharacter threories}, arXiv:1605.08971


 \bibitem{Ben}
 Benedett~C., Combinatorial Hopf algebra of supercharacters of tipe $D$, Journal of algebraic Combinatorics,  2013, vol. 38, no. 4, 767-783.

 \bibitem{Number-1}
Brumbaugh~ J.L., Bulkow~M., Fleming~ P.S., Garcia~L.A., Garcia~S.R., Karaali~G., Michal~ M., Turner~A.P., Suh~H. \emph{Supercharacters, exponentioal sums and the uncertainty principle}, Journal of Number theory, 2014, vol. 144,  151-175.

\bibitem{Burkett}
Burkett Sh., Lamar J., Lewis M.L., Wynn C.,  \emph{Groups with exactly two supercharacter theories}, arXiv:1506.00015

\bibitem{DI}
Diaconis~ P., Isaacs~I.M., \emph{Supercharacters and  superclasses for
algebra groups}, Trans.Amer.Math.Soc., 2008, vol. 360, 2359-2392.

\bibitem{Number-2}
Fowler~C.F., Garcia~S.R., Karaali~G., \emph{Ramanujan sums as supercharacters}, The Ramanujan Journal, 2014, vol.32,  205-241.

\bibitem{GR}
Grinberg~D., Reiner~V., \emph{Hopf algebras in combinatorics}, preprint arXiv:1409.8356


\bibitem{H}
Hendrickson A.O.F., {\it Supercharacter theory costructions corresponding to Schur ring products}, Comm. Algebra, 2012,
vol. 40, no.12, 4420-4438.

\bibitem{IsKar}
Isaacs~I.M.,  Karagueuzian~D., \emph{Involution and characters of upper triangular matrix groups},  Math. of Computation, 2005,  74:2027-2033.



\bibitem{Ka}
Kazhdan D.,  \emph{Proof of Springer hypothesis}, Israil J.Math., 1977, 28:272-284.

\bibitem{Kir1}
Kirillov A.A.,  Lectures on the orbit method, Graduate
Studies in Math., vol. 64, 2004.

\bibitem{Kir2}
Kirillov A.A.,  \emph{Variations on the Triangular Theme}, Amer. Math. Soc. Transl., 1995,  169:43-73

\bibitem{CR}
Curtis C.W., Reiner I., Representation theory of finite groups and associative algebras, 1962, New York, London: Interscience piblishers.

\bibitem{MT}
Marberg~E., Theim~N.,\emph{Superinduction for pattern groups}, J. Algebra, 2009, vol. 321,  3681-3703.
\bibitem{M2}
Marberg E., \emph{Combinatorial methods of character enumeration for the unitriangular group}, J. Algebra, 2011, vol. 345, 295-323
\bibitem{Otto}
Otto B., \emph{Constituents 0f supercharacters and Kirillov functions}, Archiv der Mathematik, 2010, vol.94, 319-326.

\bibitem{P1}
 Panov A.N. \emph{Supercharacter theory for groups of invertible elements of reduced algebras}, 2016, St. Petersburg Math. J. 2016, vol. 27, 1035-1047.
\bibitem{P2}
Panov~A.N., \emph{Supercharacters for the finite groups of  triangular type}, 	arXiv:1508.05767
\bibitem{P3}
Panov~A.N., \emph{Restriction and induction for supercharacters of finite groups of triangular type}, arXiv:1610.04846
\bibitem{P4}
Panov~A.N., \emph{Invariants of the coadjoint action on the basic varieties of the unitriangular group}, Transformation groups,
vol.20, no. 1, 2015, 229-246.

\bibitem{Pi}
Pierce  R.S., Associative algebras, Springer-Verlag, New York, 1982.

\bibitem{Serr}
Serr J.-P., Linear representations of finite groups, 1977,  New-York: Springer-Verlag.


\bibitem{ThiemB}
Thiem~N., \emph{Branching rules in the sing of super functions of unipotent upper-triangular matrices}, Journal of Algebraic combinatorics, 2010, vol. 31, 267-298.

\bibitem{ThiemRest}
Thiem~N., Venkateswaran~V.,\emph{Restricting supercharacters of the finite unipotent upertriangular matrices}, Electronic Journal of Combinatorics, 2009, vol. 16, paper number 23



\bibitem{Zel} Zelevinsky~A., Representations of classical groups. A Hoph algebra approach, Lecture Notes in Math. no. 869,  Springer-Verlag, Berlin-Heidelberg-New York, 1981


\bibitem{Yan}
Ning Yan,  Representation Theory of finite unipotent linear groups, Ph.D. Thesis, Department of mathematics, University of  Pennsylvania, 2001 (see also arXiv: 1004.2674).









\end{thebibliography}
\end{document}